\documentclass[12pt]{article}

\usepackage[usenames,dvipsnames]{color}

\def\red#1{{\textcolor{Red}{#1}}}      

\usepackage{amsmath, amsthm, amscd, amsfonts, amssymb}

\newtheorem{theorem}{Theorem}[section]
\newtheorem{proposition}[theorem]{Proposition}

\newtheorem{corollary}[theorem]{Corollary}
\newtheorem{lemma}[theorem]{Lemma}
\newtheorem{definition}[theorem]{Definition}

\newtheorem{remark}[theorem]{Remark}

\def\rrd{{\mathbb{R}^d}}

\def\calf{{\mathcal{F}}}
\def\calm{{\mathcal{M}}}
\def\calo{{\mathcal{O}}}

\def\cald{{\mathcal{D}}}
\def\calx{{\mathcal{X}}}

\def\call{{\mathcal{L}}}

\def\calp{{\mathcal{P}}}

\def\vsp{\vspace*{1,5mm}\\ }
\def\bk{\bigskip }
\def\mk{\medskip }
\def\sk{\smallskip }
\def\n{\noindent }
\def\dd{\displaystyle}

\def\barr{\begin{array}}
\def\earr{\end{array}}

\def\bit{\begin{itemize}}
\def\eit{\end{itemize}}

\def\FP{Fokker--Planck}

\def\1{^{-1}}
\def\one{\mbox{$1\!\!\,\rule{0,2mm}{3,1mm}\,$}}

\def\rr{{\mathbb{R}}}
\def\nn{{\mathbb{N}}}
\def\9{{\infty}}
\def\lbb{{\lambda}}
\def\wt{\widetilde}
\def\ov{\overline}
\def\vf{{\varphi}}
\def\oo{{\omega}}
\def\ooo{{\Omega}}
\def\pp{{\partial}}
\def\vp{{\varepsilon}}

\def\ff{\forall }
\def\({\left(}
\def\){\right)}
\def\<{\left<}
\def\>{\right>}

\title{Uniqueness for nonlinear Fokker--Planck equations and for McKean--Vlasov SDEs:\\ The degenerate case}
\author{Viorel Barbu\thanks{Al.I. Cuza University and Octav Mayer Institute of Mathematics of  Romanian Academy, Ia\c si, Romania.  Email: vbarbu41@gmail.com}\and Michael R\"ockner\thanks{Fakult\"at f\"ur Mathematik, Universit\"at Bielefeld,  D-33501 Bielefeld, Germany.  Email: roeckner@math.uni-bielefeld.de}}
\date{}

\begin{document}
\maketitle
\begin{abstract}
\n This work is concerned with the existence and uniqueness of gene\-ra\-lized (mild or distributional) solutions to (possibly degenerate)  Fokker--Planck equations $\rho_t-\Delta\beta(\rho)+{\rm div}(Db(\rho)\rho)=0$ in $(0,\9)\times\rr^d,$ $\rho(0,x) \equiv \rho_0(x)$. Under suitable assumptions on  $\beta:\rr\to\rr,\,b:\rr\to\rr$   and $D:\rr^d\to\rr^d$,  \mbox{$d\ge1$,} this equation generates a unique flow $\rho(t)=S(t)\rho_0:[0,\9)\to L^1(\rr^d)$   as a mild solution in the sense of  nonlinear semigroup theory.   This flow is also unique in  the class of $L^\9((0,T)\times\rrd)\cap L^\9((0,T);H^{-1}),$ $\ff T>0$,    Schwartz distributional solutions on $(0,\9)\times\rr^d$. Moreover, for $\rho_0\in L^1(\rrd)\cap H\1(\rrd)$, $t\to S(t)\rho_0$ is differentiable from the right on $[0,\9)$ in $H\1(\rrd)$-norm. As a main application, the weak uniqueness of the corres\-pon\-ding McKean--Vlasov SDEs is  proven.\sk
\end{abstract}
\section{Introduction}\label{s1}
We  shall treat here the nonlinear \FP\ equation (NFPE)
\begin{equation}\label{e1.1}
\hspace*{-3mm}\barr{l}
\rho_t(t,x){-}\Delta\beta(\rho(t,x)){+}{\rm div}(D(x)b(\rho(t,x))\rho(t,x))=0,\\\hfill (t,x)\in(0,\9)\!\times\!\rr^d,\\
\rho(0,x)=\rho_0(x),\ x\in\rr^d,
\earr\end{equation}where $\beta:\rr\to\rr$ is monotonically nondecreasing and \mbox{$ D:\rr^d\to\rr^d$, $b:\rr\to\rr$} are given functions to be made precise below in Hypotheses (i)--(iv).

The Cauchy problem \eqref{e1.1} with the conditions
\begin{eqnarray}
&\rho(t,x)\ge0,\ \ff\,t\in[0,\9)\mbox{ and a.e. }x\in\rr^d,\label{e1.2}\\
&\dd\int_{\rr^d}\rho(t,x)dx=\int_{\rr^d}\rho(0,x)dx=1,\ \ff\,t\ge0,\label{e1.3}
\end{eqnarray}is relevant in statistical mechanics (see, e.g., \cite{23a}, \cite{9a}), mean field game theory (\cite{Lions}, \cite{19a}), as well as in stochastic analysis, where it is used to reproduce the microscopic dynamics of the solution $X(t)$  to the McKean--Vlasov stochastic differential equation of Nemytskii type (also called a {\it distribution density dependent or singular McKean--Vlasov equation})
\begin{equation}\label{e1.4}
\barr{l}
dX(t)=D(X(t))b(\rho(t,X(t)))dt+\sqrt{\dd\frac{2\beta(\rho(t,X(t))}{\rho(t,X(t))}}\, dW(t),\\
X(0)=X_0,\earr
\end{equation}
by the macroscopic dynamics of its time marginal law.

In fact, if $\rho:[0,\9)\to L^1(\rr^d)$ is a distributional solution to \eqref{e1.1},  then \eqref{e1.4} has a probabilistically weak solution $X$ on some probability space $(\ooo,\calf,\mathbb{P})$ with normal filtration $(\calf_t)_{t\ge0}$ and $(\calf_t)$-Brownian motion $W$ with values in $\rr^d$ such that  $\rho(t,x)dx=\mathbb{P}\circ(X(t))\1(dx),$ $\rho_0(x)dx= \mathbb{P}\circ(X_0)\1(dx)$. (See \cite{2a}--\cite{7b}.) 
In these papers, we have developed an approach to McKean--Vlasov SDEs (see \cite{21a}, \cite{28a}) starting from the corresponding nonlinear \FP\ equations and based on the superposition principle (see \cite{10}). Our approach is, in particular, taylored to apply to Nemytskii type McKean--Vlasov SDEs as above, which -- to the best of our knowledge -- first appeared in \mbox{McKean's} fundamental paper \cite{21a}. Since then, a large literature has emerged.   In this context, we refer to the references in \cite{2a}--\cite{7b} and, in particular, to \cite{20aaa}, \cite{27a} and the recent monograph \cite{17a}. 

Our hypotheses on $\beta,b$ and $D$  are the following:
\begin{itemize}
	\item[\rm(i)] $\beta\in C^2(\rr),\ \beta'(r)>0,\ \ff\,r\ne0,\ \beta(0)=0$ and, for $\alpha_1>0$,
\begin{equation}\label{e1.5}
|\beta(r)|\le\alpha_1|r|,\ \ff\,r\in\rr.\end{equation}
	\item[\rm(ii)] $D\in L^\9(\rrd;\rrd),{\rm div}\,D\in L^{m}_{\rm loc}(\rrd),$ $({\rm div}\,D)^-\in L^\9$, where $m>\frac d2$ if~$d\ge2$, $m=1$ if $d=1.$
	\item[\rm(iii)] $b\in C^1(\rr)\cap C_b(\rr), $ $b(r)\ge0,\ \ff r\in\rr.$
	\item[\rm(iv)] $|b^*(r)-b^*(\bar r)|\le\alpha_2|\beta(r)-\beta(\bar r)|,\ \ff\,r,\bar r\in\rr,$\\ where $b^*(r)\equiv b(r)r$ and $\alpha_2>0.$
	\end{itemize} 
In statistical mechanics, \FP\ equations of the form \eqref{e1.1} arise for instance in the description of the Bose--Einstein statistic, where\break \mbox{$\beta(r)\equiv a\ln(1+r)$} or in the Plastino--Plastino model of thermostatics when $\beta(r)\equiv ar^\alpha$, $0<\alpha\le1$ which is the so called "power law" diffusion (see, e.g., \cite{23a}).  Furthermore, obviously, $b (r):=\frac{\beta(r)}{r}$, $r\in\rr$, satisfies (iv), which constitutes a very in\-te\-resting case for \eqref{e1.4}, because then the diffusion coef\-fi\-cient of the Kolmogorov operator corresponding to \eqref{e1.4} is equal to the "strength" of the drift. Since the diffusion coefficients determine the appropriate (sub)Riemannian metric for \eqref{e1.4}, this gives a geometrically very interesting interpretation of~\eqref{e1.4}. This example is also interesting in Physics, if $D=-\nabla V$, for some weakly differentiable $V:\rrd\to[0,\9)$. We recall that once we have found a solution $\rho$ to \eqref{e1.1} (e.g. by applying Theorem \ref{t2.1} below), then the stochastic dynamics we get from solving \eqref{e1.4} for this $\rho$ (e.g. by applying the main result from \cite{2}) has the (linearized) Kolmogorov operator
$$L\vf(t,x):=\frac{\beta(\rho(t,x))}{\rho(t,x)}\,\Delta\vf(x)-\frac{\beta(\rho(t,x))}{\rho(t,x)}\,\nabla V(x)\cdot\nabla\vf(x),\ t\ge0,\ \vf\in C^\9_0(\rrd),$$which can be considered as the Hamiltonian of the system. Then, one can look at the "ground state transform" of $L$ given by
$$L\to e^VL(e^{-V})=:L_V,$$and a simple calculation shows that
$$L_V\Psi=\frac{\beta(\rho)}{\rho}\,\Delta\Psi-\(\frac{\beta(\rho)}{\rho}\,\Delta V\)\Psi,\ \Psi\in e^V\cdot C^\9_0(\rrd).$$So, clearly, the multiplication operator given by $\frac{\beta(\rho)}\rho\,\Delta V=b(\rho)\Delta V$ is a time-dependent potential for this transformed (linearized) Kolmogorov operator. This shows the r\^ole of $b$ and $D=-\nabla V$ in a physical context, in this case.

In general, for $\rho_0\in L^1$, NFPE \eqref{e1.1} does not have a classical (strong) solution and the best one can expect is a generalized solution in the sense of the next definition.  

\begin{definition}\label{d1}\rm A function $\rho:[0,\9)\to L^1(\rrd)$ is said to be a {\it mild solution} to \eqref{e1.1} if $\rho\in C([0,\9);L^1(\rrd))$ and we have
	\begin{eqnarray}
\rho(t)=\dd\lim_{h\to0}\rho_h(t)\mbox{ in }L^1(\rrd),\ \ff\,t\in[0,\9),\label{e1.6}\end{eqnarray}uniformly on compacts in $[0,\9)$, 
where $\rho_h:[0,T]\to L^1(\rrd)$ is the step function,
\begin{eqnarray}
&\rho_h(t)=\rho^j_h,\ \ff\,t\in [jh,(j+1)h),\ j=0,1,...,N=\mbox{$\left[\frac Th\right]$},\label{e1.7}\\[1mm]
&\rho^{j+1}_h-h\Delta\beta(\rho^{j+1}_h)+h\,{\rm div}(Db(\rho^{j+1}_h)\rho^{j+1}_h)=\rho^{j}_h\mbox{ in }\cald'(\rrd),\label{e1.8}\\[1mm]
&\rho^j_h\in L^1(\rrd),\ \ff\,j=0,...,N;\ \rho^0_h=\rho_0.\label{e1.9}\end{eqnarray}
\end{definition}
We note that \eqref{e1.7}--\eqref{e1.9} can be equivalently written as
\begin{equation}\label{e1.10}
\barr{l}
\dd\frac1h\,(\rho_h(t))-\rho_h(t-h)-\Delta\beta(\rho_h(t))+{\rm div}(Db(\rho_h(t))\rho_h(t))=0,\\\hfill \ff\,t\in[h,T),\\
\rho_h(t)=\rho_0,\ \ff\,t\in[0,h),\earr\end{equation}where \eqref{e1.10} is meant in the sense of Schwartz distributions. If we denote by $A_0:L^1(\rrd)\to L^1(\rrd)$ the operator
\begin{equation}\label{e1.10a}
\barr{rcl}
A_0(u)&=&-\Delta\beta(u)+{\rm div}(Db(u)u),\ \ff\,u\in D(A_0),\vsp
D(A_0)&=&\left\{u\in L^1(\rrd);-\Delta\beta(u)+{\rm div}(Db(u)u)\in L^1(\rrd)\right\},\earr\end{equation}
where $\Delta$ and div are taken in the sense of Schwartz distributions on $\rrd$, then  the system \eqref{e1.8}--\eqref{e1.9} is equivalent to
\begin{equation}\label{e1.11}
\barr{l}
\rho^{j+1}_h+hA_0(\rho^{j+1}_h)=\rho^j_h,\ j=0,1,...,N,\vsp
\rho^0_h=\rho_0.\earr\end{equation}
In other terms, this means (see, e.g., \cite{1}, p.~129) that $\rho$ is a {\it mild solution} to the Cauchy problem

\begin{equation}\label{e1.12}
\barr{l}
\dd\frac{d\rho}{dt}+A_0(\rho)=0,\ \ff\,t\ge0,\vsp
\rho(0)=\rho_0.\earr\end{equation}
By the general existence theory for the nonlinear Cauchy pro\-blem  in Banach spaces (the Crandall \& Liggett existence theorem), for each $\rho_0\in\ov{D(A_0)}$, pro\-blem \eqref{e1.12} has a unique mild solution $\rho\in C([0,\9);L^1(\rrd))$, if $A_0$ is $m$-accretive in $L^1(\rrd)$, that is,
\begin{eqnarray}
&R(I+\lbb A_0)=L^1(\rrd),\ \ff\,\lbb>0,\label{e1.13}\\[2mm]
&\|(I+\lbb A_0)\1u-(I+\lbb A_0)\1v\|_{L^1(\rrd)}\le\|u-v\|_{L^1(\rrd)},\qquad\label{e1.14}\\
&\ \hspace*{55mm}\ \ff u,v\in L^1(\rr^d). \nonumber
\end{eqnarray}
Moreover, in this case (see, e.g., \cite{1}, p.~154), the solution $\rho=\rho(t,\rho_0)\equiv S(t)\rho_0$ defines a semigroup of contractions on $\ov{D(A_0)}$, that is, $\rho(t,\rho_0)\in\ov{D(A_0)}$, $\ff\,t\ge0,$ and
\begin{eqnarray}
&\rho(t+s,\rho_0)=\rho(t,\rho(s,\rho_0)),\ \ff\,t,s,\ge0,\ \rho_0\in\ov{D(A_0)},\label{e1.15}\\[2mm]
&\|\rho(t,\rho_0)-\rho(t,\bar\rho_0)\| _{L^1(\rrd)}\le\|\rho_0-\bar\rho_0\|_{L^1(\rrd)},\ \ff\,t\ge0,\ \rho_0,\bar\rho_0\in\ov{D(A_0)},\ \ \qquad\label{e1.16}
\end{eqnarray}
and, as seen by \eqref{e1.6}--\eqref{e1.9}, $\rho$    also denoted $e^{-tA_0}$, is given by the exponential formula
\begin{eqnarray}
&\rho(t,\rho_0)=\dd\lim_{n\to\9}\(I+\dd\frac tn\,A_0\)^{-n}\rho_0\ \ \mbox{ in }L^1(\rrd),\ \ff\,t\ge0,\qquad\label{e1.17}
\end{eqnarray}uniformly on compacts in $[0,\9)$. 	In the works \cite{2}--\cite{5}, under even weaker hypotheses than (i)--(iii), the range condition \eqref{e1.13} and the existence of a family $\{J_\lbb\}$ of nonlinear contractions in $L^1(\rrd)$ with $J_\lbb(L^1(\rr^d))\subset D(A_0)$ was proven, which satisfy the resolvent equation
$$J_{\lbb_2}(f)=J_{\lbb_1}\(\frac{\lbb_1}{\lbb_2}\,f+\(1-\frac{\lbb_1}{\lbb_2}\)J_{\lbb_2}(f)\),\ \ff\,\lbb_1,\lbb_2>0,\ f\in L^1,$$ and $J_\lbb(f)\in(I+\lbb A_0)\1f$, $\ff f\in L^1(\rrd)$, $\lbb>0$.

Then, the   operator $A: {D(A)}\subset D(A_0)\subset L^1(\rrd)\to L^1(\rrd)$,  defined~by
\begin{equation}\label{e1.18}
\barr{rcl}
A(u)&\!\!=\!\!&A_0(u),\ \ u\,=\,J_\lbb(f),\ \  f\in L^1(\rrd),\vsp
D(A)&\!\!=\!\!&J_\lbb(L^1(\rrd)),\earr\end{equation}
is independent of $\lbb$,   $m$-accretive and
\begin{equation}\label{e1.19}
(I+\lbb A)\1(f)=J_\lbb(f)\in(I+\lbb A_0)\1f,\ \ff f\in L^1(\rrd),\ \lbb>0.\end{equation}
This means that, for $\rho_0\in\ov{D(A)}$, the mild solution $\rho$ to the Cauchy problem
\begin{equation}\label{e1.20}
\barr{l}
\dd\frac{d\rho}{dt}+A(\rho)=0,\ \ t\ge0,\vsp
\rho(0)=\rho_0,\earr\end{equation}is just a mild solution to \eqref{e1.1} in the sense of Definition \ref{d1}.  However, this does not imply the uniqueness of the mild solution to \eqref{e1.1}, because the family $\{J_\lbb\}$,  and so the operator $A$, which is defined by $\{J_\lbb\}$,  are not unique. In fact, $J_\lbb(f)$ is given~by
\begin{equation}\label{e1.21a}
J_\lbb(f)=\lim_{\vp\to0} y_\vp\mbox{ in }L^1(\rrd),\end{equation}where $y_\vp\in L^1(\rrd)\cap H^2(\rrd)$ is the solution to an approximating equation of the form
\begin{equation}\label{e1.21}
y_\vp-\lbb(\vp I+\Delta)\beta(y_\vp)+\lbb\,{\rm div}(Db_\vp(y_\vp)y_\vp)=f\mbox{ in }\cald'(\rrd),\end{equation}where $b_\vp$ is a smooth approximation of $b$. (Other approximating equations of the form \eqref{e1.21} could be considered as well.)  However, since the limit \eqref{e1.21a} depends on the sequence $\{y_\vp\}$ and hence might be not unique, one obtains in this way a family of mappings $\{J_\lbb\}_{\lbb>0}$ and each one defines via \eqref{e1.18} an $m$-accretive operator $A$. So the uniqueness  is not in the class of all mild solutions $\rho$ to \eqref{e1.20} and for uniqueness one should request a further condition to restrict the class of mild solutions. Such a situation is encountered for the conservation law equation as well, i.e., for $\beta\equiv0$, where the corresponding mild solution is unique in the narrow class of Kruzkov entropic solutions \cite{11a}.

The aim of this paper is twofold. The first is to prove that, if $1\le d$, then {\it under the additional  Hypothesis {\rm(iv)} the operator $A_0$ is itself $m$-accretive} (that is, \mbox{$I+\lbb A_0$} is invertible for all $\lbb>0$, and also \eqref{e1.14} holds) and so  {\it the Crandall \& Liggett existence theorem is applicable to the Cauchy problem \eqref{e1.12} to derive not only the existence, but also the uniqueness of a mild solution $\rho$ to \eqref{e1.1}} (Theo\-rem~\ref{t2.1}).
The second is to prove the uniqueness of solutions $\rho$ to \eqref{e1.1} in the class of $L^1$-valued  distributional solutions $\rho$.   This means $\beta(\rho)\in L^1_{\rm loc}((0,\9)\times\rrd)$ and 
\begin{equation}\label{e1.22}
\hspace*{-5mm}\barr{l}
\dd\int^\9_0\!\!\!\int_\rrd\!(\rho\vf_t+\beta(\rho)\Delta
\vf+b(\rho)\rho(D\cdot\nabla\vf))dt\,dx
+\int_{\rr^d}\vf(0,x)\rho_0(dx)=0,\\
\hfill\ff\vf\in C^\9_0([0,\9)\times\rrd),\earr\end{equation}where $\rho_0$  is a signed Radon measure on $\rr^d$ of bounded variation.  
As seen later on in Theorem \ref{t2.1}, under the above assumptions the {\it mild solution $\rho$ is also a distributional solution to \eqref{e1.1}}. However, the uniqueness of distributional  solutions to \eqref{e1.1} is still an open problem  under the general Hypotheses (i)--(iii). For the porous media equation (that is, \mbox{$D\equiv0$}) such a uniqueness result was established by H.~Brezis and M.G.~Crandall \cite{8} and M.~Pierre \cite{9} in the class   $L^\9((0,\9)\times\rrd)$ for distributional solutions and a similar result was established in \cite{5a}, \cite{7b} for NFPE \eqref{e1.1} if $\beta$ is strictly monotone. (See, also, \cite{7} for $d=1$.)  Here, one proves such a result if one merely assumes that $\beta$, $D$ and $b$ sa\-tisfy Hypotheses (j)--(jv) in Section \ref{s3}  (Theorem~\ref{t3.1} and Corollary \ref{c3.2''}). In particular,   degenerate cases where $\{\beta'=0\}\ne\emptyset$ are covered and our results generalize \cite{8} to the case $D\not\equiv0.$

In the limit case $\beta\equiv0$, $Db(r)\equiv\{a_i(r)\}^d_{i=1}\equiv a(r)$, equation \eqref{e1.1} reduces~to
\begin{equation}\label{e3.11}
	\barr{l}
	\rho_t+{\rm div}(a(\rho))=0,\ \cald'((0,\9)\times\rrd),\vsp
	\rho(0)=\rho_0,\earr
	\end{equation}
which has a unique Kru\v zkov's solution $\rho(t)=S_0(t)\rho_0$ defined as
\begin{equation}
\label{e3.12}
|\rho-k|_t+{\rm div}((a(\rho)-a(k)){\rm sign}(\rho-k))\le0\mbox{ in }\cald'((0,\9)\times\rrd),\ k\in\mathbb{Z}^+,
\end{equation}
$\rho\in C([0,\9);L^1),\ \rho(0)=\rho_0.$

So, contrary to the situation encountered in Theorem \ref{t3.1}, the uniqueness for \eqref{e3.11} is not in the class of distributional solutions but in that of entropic solutions and this fact emphasizes the role of the diffusion term $\Delta\beta(u)$ in distributional uniqueness.

The notion of entropy solution was extended by Carillo \cite{18a}, Chen and Pertham\'e \cite{20ay} to NFPE \eqref{e1.1} and to more general degenerate parabolic--hyperbolic equations. It should be mentioned, however, that in general these two concepts, that is, mild and entropy solutions are different if $u$ is not in $L^\9(0,\9;L^1(\rrd))$ and so the uniqueness of a mild solution is not covered by the above mentioned results. Of course, this happens also for the more general class of distributional solutions to \eqref{e1.1}.

We also stress here that the class of distributional solutions, i.e. the solutions to \eqref{e1.22} is a much wider class of solutions than the class of entropy solutions, as e.g. the latter are weakly differentiable in space by definition, while the solutions to \eqref{e1.22} are not in general. Therefore, the uniqueness results for distributional solutions in Section \ref{s3} are considerably stronger than results on uniqueness results for  entropy solutions. Furthermore, the latter would not be sufficient for the applications to obtain uniqueness of weak solutions to the McKean--Vlasov equation \eqref{e1.4}, since (in particular, in the degenerate case) there is no reason to expect that the time-marginal laws of a solution to \eqref{e1.4} are always weakly differentiable in space. Our  uniqueness results on weak solutions to \eqref{e1.4} in Section \ref{s.4} are, therefore, based on our uniqueness results on distributional solutions. More precisely, by virtue of the equivalence of equations \eqref{e1.1} and \eqref{e1.4}, the above results have implications on the uniqueness of probabilistically weak solutions to the McKean-Vlasov equation \eqref{e1.4}. The key additional analytical result is to  prove "linearized uniqueness" for \eqref{e1.22} (see Theorem \ref{t4.1} and Corollary~\ref{c4.2}), which together with Corollary \ref{c3.2''} imply weak uniqueness for SDE \eqref{e1.4} (see Theorem \ref{t4.3} and Remark \ref{r4.5a} below). 

Concerning other work on weak uniqueness of McKean--Vlasov SDEs, we refer to \cite{16a}, \cite{20a}, \cite{20aa}, \cite{21aa} and \cite{21aaa}. However, none of them covers our results in Theorem \ref{t4.3} below, because either the diffusion coefficients do not depend on the time marginal laws or other additional assumptions (as, e.g., linear growth) are assumed. In addition, all these papers assume nondegeneracy conditions on the noise, whereas our Theorem \ref{t4.3} allows degenerate diffusivity.

As regards the literature on generalized (mild) solutions to nonlinear \FP\ equations via the nonlinear semigroup theory, we refer to   \cite{2}--\cite{4}, \cite{5}, \cite{8c}   and the related papers   \cite{12b}, \cite{13b}.)

The time-dependent case was treated in \cite{8b} under time-regularity hypotheses on the diffusion and drift coefficients invoking the general existence theory for the Cauchy problem in a Banach space with time-dependent accretive operators.

\bk\noindent{\bf Notation.} $L^p(\rrd),\ 1\le p\le\9$ (denoted $L^p$) is the space of all Lebesgue measurable and $p$-integrable functions on $\rr^d$, with the standard norm \mbox{$|\cdot|_p$.} $(\cdot,\cdot)_2$ denotes the inner product in $L^2$. By $L^p_{\rm loc}$ we denote the correspon\-ding local space.  For any open set $\calo\subset\rrd$ let    $W^{k,p}(\calo)$, $k\ge1$, denote the standard Sobolev space on  $\calo$ and by $W^{k,p}_{\rm loc}(\calo)$ the corresponding local space. We set $W^{1,2}(\calo)=H^1(\calo)$,
$W^{2,2}(\calo)=H^2(\calo)$, $H^1_0(\calo)=\{u\in H^1(\calo),$ $u=0\mbox{ on }\pp\calo\}$, where $\pp\calo$ is the boundary of $\calo$. By $H\1(\calo)$ we denote the dual space of $H^1_0(\calo)$ (of $H^1(\rrd)$, respectively, if \mbox{$\calo=\rrd$).} In the following, we shall simply write $H^1=H^1(\rrd)$, $H\1=H\1(\rrd)$. $C^\9_0(\calo)$ is the space of infinitely differentiable real-valued functions with compact support in $\calo$ and $\cald'(\calo)$ is the dual of $C^\9_0(\calo)$, that is, the space of Schwartz distributions on $\calo$.  Denote by $C^k(\rr)$ the space of all  continuously differentiable real-valued functions on $\rr$ up to order $k$, by $C_b(\rr)$ the space of continuous and bounded real-valued functions on $\rr$, and by ${\rm Lip}(\rr)$ the space of real-valued Lipschitz functions on $\rr$ with the norm denoted by $\|\cdot\|_{\rm Lip}$.  $C([0,\9);L^1(\rrd))$ is the space of continuous functions $y:[0,\9)\to L^1(\rrd).$ Denote also by $C^\9_0([0,\9)\times\rr^d)$ the space of all $\vf\in C^\9([0,\9)\times\rr^d)$ such that \mbox{$support\, \vf\subset K$,} where $K$ is compact in $[0,\9)\times\rr^d$. We shall also use the following notations: 
$$\barr{c}
\beta'(r)\equiv\dd\frac d{dr}\,\beta(r),\ b'(r)\equiv\dd\frac d{dr}\,b(r),
\vsp
y_t=\dd\frac\pp{\pp t}\,y,\ \Delta y=\!\!\dd\sum^d_{i=1}\frac{\pp^2}{\pp x^2_i}\,y,\ \nabla y=\!\left\{\dd\frac{\pp y}{\pp x_i}\right\}^d_{i=1}\!\!\!,\
{\rm div}\,u=\!\sum^d_{i=1} \dd\frac{\pp u_i}{\pp x_i}  ,\ u=\{u_i\}^d_{i=1}\earr$$We also denote by $\calp$ the set of all probability densities $\rho$ on $\rr^d$, that is,
$$\calp=\left\{y\in L^1(\rrd);\ y\ge0,\mbox{ a.e. in }\rrd,\ \int_\rrd y(x)dx=1\right\}.$$
Denote by $\calm(\rr^d)$ the space of all signed Radon measures on $\rr^d$of bounded variation.  A sequence $\{\mu_n\}\subset\calm(\rr^d)$ is said to be converging to $\mu$ in $\sigma(\calm(\rr^d),C_b(\rr^d))$ topology~if
\begin{equation}
\label{1.26a}
\lim_{n\to\9}\int_{\rr^d}\psi\,d\mu_n=\int_{\rr^d}\psi\,d\mu,\ \ \ff\psi\in C_b(\rr^d).
\end{equation}
The function $t\to\mu(t)\in\calm(\rr^d)$ is said to be {\it narrowly continuous} on $[0,\9)$ if, for every $\psi\in C_b(\rr^d)$, the function  $t\to\int_{\rr^d}\psi\,d\mu(t)$ is continuous on on~$[0,\9)$.

 \section{The existence and   uniqueness of a mild\\  solution to NFPE}\label{s2}
 \setcounter{equation}{0}

Theorem \ref{t2.1} is one of the main results of this work.

 \begin{theorem}\label{t2.1} Assume that $1\le d.$ Then, under Hypotheses {\rm(i)--(iv)}, for each $\rho_0\in L^1$ there is a unique mild solution $\rho=\rho(t,\rho_0)$ to equation \eqref{e1.1} which is also a distributional solution. Moreover, \eqref{e1.15}--\eqref{e1.17} hold and, if $\rho_0\in\calp$, then $\rho(t)\in\calp$, $\ff\,t\ge0$, that is, \eqref{e1.2}--\eqref{e1.3} hold. 
 	
 	Finally, if $\rho_0\in L^1\cap L^\9$, then $\rho\in L^\9((0,T)\times\rr^d)$ for all $T>0$.
 \end{theorem}
By definition, the continuous semigroup $t\to\rho=\rho(t,\rho_0)$ given by Theo\-rem \ref{t2.1} is the {\it nonlinear \FP\ flow} associated with the nonlinear diffusion $\beta=\beta(\rho)$ and the
drift term $\rho\to Db(\rho)\rho$. This means that the operator $A_0$ defined by \eqref{e1.10a}
is the infinitesimal generator of the continuous semigroup of contractions $S(t)\rho_0=\rho(t,\rho_0)$, that is, of the \FP\ flow $S(t)$. Theorem \ref{t2.1} can be rephrased in terms of   semigroup theory as follows: {\it Under Hypotheses {\rm(i)--(iv)}, the operator $A_0$ generates in a weak-mild sense a con\-ti\-nuous semigroup of contractions in $L^1(\rrd)$ which leaves invariant the set $\calp$ of all probability densities.} 
Moreover, by Theorem 5.1 in \cite{4} it follows that, if $d=3$ and $\beta'(r)\ge a>0$, then the flow $S(t)$ has a smoothing effect on initial data and extends to all $\rho_0\in\mathcal{M}_b$.   

It should be emphasized that such a generation result for the \FP\ flow $S(t)$ is specific to  the space $L^1(\rrd)$, because only in this space the operator $A_0$ is accretive. Moreover, the semigroup $S(t)$ is not $t$-differentiable in this space and so the solution $\rho$ exists in the above mild sense only.  However, for  $\rho_0\in L^1\cap L^2$,   $\rho$   is a strong solution to \eqref{e1.1} in $H\1(\rrd)$. In~fact, we~have
\begin{theorem}
\label{t2.2}
Let Hypotheses {\rm (i)--(iv)} hold and assume in addition that
\begin{itemize}
	\item[\rm(v)] $\beta(r)r\ge\alpha_0|r|^2,\ \ \ff r\in\rr,$ where $\alpha_0>0.$
\end{itemize}Assume that $\rho_0\in L^1\cap L^2$.  Then, every mild solution $\rho$ to~\eqref{e1.1} satisfies
\begin{equation}\label{e2.1b}
	\beta(\rho)\in L^2(0,T;H^1(\rrd))\cap L^\9(0,T;L^2),\ \ \ff\,T>0.\end{equation}
Moreover, $\rho:[0,T]\to H\1(\rrd)$ is absolutely continuous and
\begin{eqnarray}
	&\dd\frac{d\rho}{dt}\in L^2(0,T;H\1(\rrd)),\label{e2.2b}\\
	&\dd\frac{d\rho}{dt}\,(t)+A_0(\rho(t))=0,\mbox{\ \ a.e. }t\in(0,T),\ \ff\,T>0.\label{e2.3b}\end{eqnarray}
\end{theorem} 
\n Theorem \ref{t2.2} amounts to saying that, a.e. $t>0$,  the semigroup $S(t)$ maps the space $L^1\cap L^2$ into $\{\rho_0\in L^1\cap L^2;\beta(\rho_0)\in H^1(\rrd)\}$ and it is a.e. $H\1(\rrd)$-valued dif\-fe\-ren\-tiable on $(0,\9)$. More will be said about this in Section \ref{s4} (Theorem~\ref{t4.1b}). 

\bk\n{\bf Proof of Theorem \ref{t2.1}.} \mk

\n As mentioned earlier, Theorem \ref{t2.1} is implied by the following key result.

\begin{proposition}\label{p2.2} Let $1\le d.$ Then, under Hypotheses {\rm(i)--(iv)}, the operator $A_0$ is $m$-accretive in $L^1$. Moreover, one has	
	\begin{eqnarray}
	&\ov{D(A_0)}=L^1,\label{e2.1}\\[1mm]
&(I+\lbb A_0)\1(\calp)\subset\calp,\ \ff\,\lbb>0.\label{e2.2}
	\end{eqnarray}
	$($Here, $\ov{D(A_0)}$ is the closure of $\ov{D(A_0)}$ in $L^1$.$)$
\end{proposition}

We shall prove Proposition \ref{p2.2} following several steps and the first one is the following uniqueness result for the stationary (resolvent) equation associated with the operator $A_0$.

\begin{lemma}
	\label{l2.3} For all $f\in L^1$ and $\lbb>0$, there is at most one solution $y\in D(A_0)$ to the equation
	\begin{equation}
	\label{e2.3}y+\lbb A_0(y)=f.
	\end{equation}
\end{lemma}
\n{\bf Proof.} We shall prove first that each solution $y\in D(A_0)$ to \eqref{e2.3} is regular. Namely,
 	\begin{eqnarray}
 \label{e2.4}\beta(y)\in W^{1,q}_{\rm loc}(\rrd),\ \ff q\in\mbox{$\left[1,\frac d{d-1}\)$},&&\mbox{if }d\ge2,\\[1mm]
\beta(y)\in W^{1,\9}_{\rm loc}(\rr^d),&&\mbox{if }d=1,\label{e2.7a}
 \end{eqnarray}and setting $\frac d{d-2}:=\9$ if $d=2$
 \begin{equation}\label{e2.7prim}
 \barr{r}
 \|\nabla\beta(y)\|_{L^q(B_R)}+\|\beta(y)\|_{L^p(B_R)}\le C_R(|f|_1+|y|_1),\ 1<p<\frac d{d-2},\vsp \ff q\in\mbox{$\left[1,\frac d{d-1}\)$},\earr
 \end{equation}for all $B_R=\{x\in\rr^d;\ |x|<\rr\}$.

 We note that \eqref{e2.7prim} follows by \eqref{e2.4} via the Sobolev--Gagliardo--Nirenberg inequality invoked below (see, e.g., \cite{9b}, p. 278).

 Consider first the case $d=1$. We have $\beta(y)\in L^1$ by Hypothesis~(i) and
 \begin{equation}\label{e2.8prim}
 (\beta(y))''-(Db^*(y))'=\frac1\lbb(y-f)\mbox{ in }\cald'(\rr).\end{equation}

 \n Since $Db^*(y)\in L^1(\rr)$, it follows that $(\beta(y))'\in L^1_{\rm loc}(\rr)$ and so \mbox{$\beta(y)\in W^{1,1}_{\rm loc}(\rr)$.} Then, by (iv), $b^*(y)\in W^{1,1}_{\rm loc}(\rr^d)$ and so, by \eqref{e2.8prim} and (ii)  we infer that $Db^*(y)\in W^{1,1}_{\rm loc}(\rr^d)$ and $(\beta(y))''\in L^1_{\rm loc}(\rr^d)$. Hence,
 $\beta(y)\in W^{1,\9}_{\rm loc}(\rr^d)$, as claimed.

Consider now the case $2\le d$. By \eqref{e2.3}, we have, for all $\vf\in C^2_0(\rrd)$,
 $$\barr{ll}
 \dd\frac1\lbb(y{-}f)\vf\!\!\!&=\!\vf\Delta\beta(y){-}\vf\,{\rm div}(Db^*(y))\vsp
 &=\!\Delta(\vf\beta(y)){-}\beta(y)\Delta\vf{-}{\rm div}(D\vf b^*(y)){-}2\nabla\beta(y)\!\cdot\!\nabla\vf{+}(D{\cdot}\nabla\vf)b^*(y)\vsp
 &=\!\Delta(\vf\beta(y)){+}\beta(y)\Delta\vf{+}(D{\cdot}\nabla\vf)b^*(y){-}{\rm div}(2\beta(y)\nabla\vf{+}D\vf b^*(y)),
  \earr $$
 and, therefore,
 	\begin{equation}
 \label{e2.5}\Delta(\vf\beta(y))=f_1+{\rm div}\,f_2\mbox{\ \ in }\cald'(\rrd),
 \end{equation}where
 \begin{equation}\label{e2.11a}
 \barr{rcl}
 f_1&=&\frac1\lbb\,(y-f)\vf-\beta(y)\Delta\vf-(D\cdot\nabla\vf)b^*(y)\vspace*{2,5mm}\\
 f_2&=&2\beta(y)\nabla\vf+D\vf b^*(y).\earr\end{equation} 
 We set $u=\vf\beta(y),\ u_\vp=u*\Psi_\vp,\ f^\vp_i=f_i*\Psi_\vp$, $i=1,2$, where $\Psi_\vp$ is a standard mollifier, that is,
  $$\Psi_\vp(x)=\frac1{\vp^d}\,\Psi\(\frac x\vp\),\ \Psi\in C^\9_0(\rrd),{\rm support}\,\Psi\subset\{x;|x|\le1\},\ \int_{\rr^d}\Psi(x)dx=1.$$Let $\calo,\calo'$ be open balls    in $\rr^d$ centered at zero such that $\ov\calo'\subset\calo$ and choose $\vf\in C^2_0(\rr^d)$ such that $\vf=1$ on $\calo'$ and  $({\rm supp}\ \vf)_\vp\subset\calo,$ $\vp\in(0,1]$, where $({\rm supp}\,\vf)_\vp$ denotes the closed $\vp$-neighbourhood of ${\rm supp}\,\vf$.
 Then, by \eqref{e2.5} we~have
 \begin{equation}\label{e2.6}
 \Delta u_\vp=f^\vp_1+{\rm div}\,f^\vp_2\mbox{\ \ in }\calo,\ \ \ u_\vp\in C^\9_0(\calo).
 \end{equation}
 Hence, by the uniqueness of the solution $u_\vp$ to \eqref{e2.6}, we have \mbox{$u_\vp=u^1_\vp+u^2_\vp$,} where $u^1_\vp,u^2_\vp\in C^\9(\calo)\cap C(\ov\calo)$ are the solutions  to the boundary value pro\-blems
 \begin{eqnarray}
 \Delta u^1_\vp=f^\vp_1\mbox{ in }\calo,&&u^1_\vp=0\mbox{ on }\pp\calo,\label{e2.7}\\[2mm]
  \Delta u^2_\vp={\rm div}\,f^\vp_2\mbox{ in }\calo,&&u^2_\vp=0\mbox{ on }\pp\calo.\label{e2.8}
 \end{eqnarray}
 By the standard existence theory for elliptic equations, we know that (see, e.g.,~\cite{8a}, Corollary 12)
   \begin{equation}\label{e2.9}
 \|u^1_\vp\|_{W^{1,q}_0(\calo)}\le C\|f^\vp_1\|_{L^1(\calo)}\le C(|y|_1+|f|_1),\ \ff\,\vp>0,\end{equation}where $1\le q<\frac d{d-1}$ and where we used Hypotheses (i)--(iii) for the last inequa\-lity. By the Sobolev--Galiardo--Nirenberg theorem (see, e.g., \cite{9b}, p.~278 and p.~281), it follows by \eqref{e2.9} that  we have 
 
 \begin{eqnarray}
 &\label{e2.9a}
 \|u^1_\vp\|_p\le C(|f|_1+|y|_1),\ \ \ff\,p\in\mbox{$\left[1,\dd\frac d{d-2}\right)$}\mbox{\ \ if }d>2,\\[1mm]
 &\|u^1_\vp\|_p\le C(|f|_1+|y|_1),\ \ \ff\,p\ge1\mbox{\ \ if }d=2.\label{e2.14a}
 \end{eqnarray}
 (In the following, we shall denote by the same symbol $C$ several positive constants independent of $\vp$ and $\|\cdot\|_p$ is the norm of $L^p(\calo)$.)

 Consider now  the solution $u^2_\vp$ to equation \eqref{e2.8}.

  If $\psi\in L^m(\calo),\ m>d$, and $\theta\in W^{2,m}(\calo)\cap W^{1,m}_0(\calo)$ is the solution to the Dirichlet problem
\begin{equation}\label{e2.16a}
-\Delta\theta=\psi\mbox{ in }\calo;\ \ \theta=0\mbox{ on }\pp\calo,\end{equation}we see by \eqref{e2.8} and by the Morrey embedding theorem (see  \cite{9b}, p.~282) that $\nabla\theta\in L^\9(\calo)$ and, therefore, by Green's formula, since $u^2_\vp=\theta=0$ on~$\pp\calo$, 
 \begin{equation}\label{e2.11aa}
 \int_\calo u^2_\vp\Delta\theta dx=-\int_\calo f^\vp_2\cdot\nabla\theta dx\le|f^\vp_2|_1\|\nabla\theta\|_\9\le C(|f|_1+|y|_1)\|\psi\|_m.\end{equation}This yields
 \begin{equation}\label{e2.23a}
 \left|\int_\calo u^2_\vp\psi dx\right|\le C(|f|_1+|y|_1)\|\psi\|_m,\ \ff\psi\in L^m(\calo).\end{equation}
Then, if   $\frac1{m'}=1-\frac1m$,   by \eqref{e2.23a} it follows by duality that
  $u^2_\vp\in L^{m'}(\calo)\subset L^q(\calo)$ for all $q\in\left[1,\frac d{d-1}\right)$ and
  $$\|u^2_\vp\|_q\le C(|f|_1+|y|_1),\ \ \ff\,\vp>0,$$and so, by \eqref{e2.9}, it follows also that $u^i_\vp\in L^{q}(\calo)$, $i=1,2$,    and
\begin{equation*}
\|u^i_\vp\|_{q}\le C(|f|_1+|y|_1),\ \ff\,q\in\mbox{$\left[1,\frac d{d-1}\right)$},\ i=1,2.\end{equation*}Hence,
\begin{equation}\label{e2.17a}
\|u_\vp\|_{q}\le C(|f|_1+|y|_1),\ \ff\,\vp>0,\   q\in\mbox{$\left[1,\frac d{d-1}\right)$}.\end{equation}
Finally, taking into account that $u_\vp=(\vf\beta(y))*\Psi_\vp$, by letting $\vp\to0$ we see by \eqref{e2.17a} and \eqref{e1.5} that
 $$\|\vf\beta(y)\|_{q}\le C(|f|_1+|y|_1),\ \ \ff\,q\in \mbox{$\left[1,\frac d{d-1}\right)$}.$$Because $\vf$ and the corresponding ball $\calo$ are arbitrary, we conclude that $y,$ $\beta(y)\in L^q_{\rm loc}(\rrd)$ and that (for a possible larger $C$, still independent of $\vp$)\newpage 
 $$\|\beta(y)\|_q\le C(|f|_1+|y|_1),\ \ff q\in\mbox{$\left[1,\frac d{d-1}\right)$}.$$  In particular, by Hypothesis (iv), this implies that  $$\|f_2\|_q\le C(|f|_1+|y|_1)$$ and, therefore,
 \begin{equation}\label{e2.11aaa}
\|f^\vp_2\|_{q}\le C(|f|_1+|y|_1),\ \ff\vp>0, q\in \mbox{$\left[1,\frac d{d-1}\right)$}. \end{equation}
Now, we shall improve the last estimate by invoking a bootstrap argument. Namely, we take in \eqref{e2.16a} $\psi\in L^\ell(\calo)$, where $\frac d2<\ell.$ This yields as above that
$$\left|\int_\calo u^2_\vp\psi\,dx\right|\le\int_\calo|f^\vp_2|\,|\nabla\theta|dx\le C\|f^\vp_2\|_q\|\nabla\theta\|_{q'}$$for all  $q\in\left[1,\frac d{d-1}\right)$ and   $q'=\frac q{q-1}>d.$ Again by the Sobolev inequality we have, for all $\ell\in\(\frac d2,d\)$ that $q':=d\ell/(d-\ell)>d$ and
\begin{equation}\label{e2.19a}
\|\nabla\theta\|_{q'}\le C\|\theta\|_{W^{2,\ell}(\calo)}\le C\|\psi\|_\ell.\end{equation}
This yields, for all $\ell\in\(\frac d2,d\)$,
\begin{equation}\label{e2.19aa}
\left|\int_\calo u^2_\vp\psi\,dx\right|\le C\|f^\vp_2\|_q\|\psi\|_\ell
\le C(|f|_1+|y|_1)\|\psi\|_\ell,\ \ff\psi\in L^\ell(\calo),\end{equation}
and, therefore, putting $\frac d{d-2}:=\9$, if $d=2$,
\begin{equation}
\label{e2.18a}\|u^2_\vp\|_r\le C(|f|_1+|y|_1),\ \ff r\in\mbox{$\left[1,\frac d{d-2}\right)$},
\end{equation}
Then, by \eqref{e2.9a}, \eqref{e2.14a},
we get
\begin{equation}\label{e2.19aaaa}
\|u_\vp\|_r\le C(|f|_1+|y|_1),\ \ff\,r\in\mbox{$\left[1,\frac d{d-2}\right)$}.\end{equation} 
Letting $\vp\to0$, this yields the bound for $\|y\|_{L^p(B_R)}$ in \eqref{e2.7prim}, since $\vf\in C^2_0(\rr^d)$ was arbitrary. Furthermore, \eqref{e2.19aaaa} implies that \eqref{e2.11aaa} is strengthened~to
\begin{equation}\label{e2.20a}
 \|f^\vp_2\|_\nu\le C(|f|_1+|y|_1),\ \ff\nu\in\mbox{$\left[1,\frac d{d-2}\right)$},\ \vp>0.\end{equation}For $\vp\to0$, this yields, since $f_2$ has compact support in $\calo$,
 \begin{equation}
 \label{e2.27prim}
 |f_2|_\nu\le C(|f|_1+|y|_1),\ \ \ff\nu\in\mbox{$\left[1,\frac d{d-2}\right)$}.
 \end{equation}Hence, for $d\in[2,3]$, we get
\begin{equation}\label{e2.20aa}
\|f^\vp_2\|_2\le C(|f|_1+|y|_1),\ \ff\vp>0,\end{equation}and so
\begin{equation*}
\|{\rm div}(f^\vp_2)\|_{H^{-1}(\calo)}\le C(|f|_1+|y|_1),\ \ff\vp>0.
\end{equation*}Then, by equation \eqref{e2.8} and, since $\frac d{d-1}\le2$, it follows that
 $u^\vp_2\in H^1_0(\calo)$ and
\begin{equation}\label{e2.28a}
\|u^2_\vp\|_{W^{1,q}(\calo)}\le\|u^2_\vp\|_{H^1_0(\calo)}
\le C(|f|_1+|y|_1),\ \ff\vp>0.\end{equation}Hence, by \eqref{e2.9}, we have
\begin{equation}\label{e2.18aa}
\|u_\vp\|_{W^{1,q}_0(\calo)}\le C(|y|_1+|f|_1),\ \ff\,\vp>0, \mbox{ if }d=2,3.\end{equation}
Letting $\vp\to0$   we get  the estimate
 \begin{equation}\label{e2.30a}
 \|u\|_{W^{1,q}(\calo)}\le C(|f|_1+|y|_1),\ \ d=2,3.\end{equation}Recalling that $u=\vf\beta(y)$, since $\vf$ and the corresponding ball $\calo$ are arbitrary, this implies \eqref{e2.4}, for $d\in[2,3]$, as claimed.

 We shall consider now the case $d\ge3.$ To this end, we come back to equation \eqref{e2.5} and note that $\vf\beta(y)=u_1+u_2$, where $u_1,u_2$ are solutions to the equations
 \begin{eqnarray}
 &&\Delta u_1=f_1\mbox{ in }\cald'(\rr^d),\label{e2.31prim}\\
 &&\Delta u_2={\rm div}(f_2)\mbox{ in }\cald'(\rr^d),\label{e2.31secund}
 \end{eqnarray}where $f_1,f_2$ are defined by \eqref{e2.11a}.

 Since $f_1\in L^1$, it follows by \cite[Lemma A.5]{9''} that $u_1$ is given by the representation formula
 $$u_1=-E*f_1\mbox{\ \ in }\rr^d,$$where $E(x)\equiv\frac1{(d-2)\oo_d|x|^{d-2}}$ is the fundamental solution to $\Delta$. 

 Hence  (see, e.g, \cite{9''}), $u_1\in M^{\frac d{d-2}}(\rr^d)\subset L^{p}_{\rm loc}(\rr^d)$, $\ff p\in\left[1,\frac d{d-2}\right)$ and $|\nabla u_1|=|\nabla E*f_1|\in M^{\frac d{d-1}}(\rr^d)\subset L^p_{\rm loc}(\rr^d)$, $\ff p\in \left[1,\frac d{d-1}\)$ with
 \begin{equation}\label{e2.31tert}
 	\|\nabla u_1\|_{L^{p}(B_R)}\le C(|y|_1+|f|_1),\ \ff R>0,\ p\in\mbox{$\left[1,\frac d{d-1}\)$}.\end{equation} (Here, $M^{\ell}$ is the Marcinkievicz space of order $\ell$.)
 	
 As regards the solution $u_2$ to equation \eqref{e2.31secund}, we note that from \eqref{e2.27prim} it follows that $f_2\in L^{p}(\rr^d)$, $\ff p\in\left[1,\frac d{d-1}\)$. Let $u_{2,\vp}=u_2*\Psi_\vp$. Then
 $$u_{2,\vp}=-\nabla E*f^\vp_2\mbox{ in }\rr^d.$$Taking into account  that $|\nabla^2 E(x)|\le C|x|^{-d},\ \ff x\ne0$, it follows by the Calderon--Zygmund theorem (see, e.g, \cite{14a} and estimate \eqref{e2.27prim}) that
 $$|\nabla u_{2,\vp}|_{p}\le C|f_2^\vp|_p\le C(|f|_1+|y|_1),\ \ff p\in\mbox{$\left[1,\frac d{d-1}\)$}$$and, after letting $\vp\to0$,   together with \eqref{e2.31tert} this yields
 $$\|\nabla u\|_{L^{p}(B_R)}\le C_R(|f|_1+|y|_1),\ \ff p\in\mbox{$\left[1,\frac d{d-1}\)$},$$and so \eqref{e2.4} and \eqref{e2.7prim}  hold  for all $d\ge2$.

 Now, let us prove the uniqueness of the solution to \eqref{e2.3}.

 If $y_1,y_2\in D(A_0)$ are two solutions,
 we have
\begin{equation}\label{e2.31a}
y_1-y_2-\lbb\Delta(\beta(y_1)-\beta(y_2))+\lbb\,{\rm div}(D(b^*(y_1)-b^*(y_2)))=0.\end{equation}

 Let $\eta\in C^2([0,\9))$ be such that
 $$\eta(r)\ge0,\ \eta(r)=1,\ \ff r\in[0,1];\   \eta(r)=0,\ \ff r\,\in[2,\9).$$
 We set
$$\vf_n(x)=\eta\(\frac{|x|^2}n\),\ \ff x\in\rrd,\ n\in\nn,$$and note that
 \begin{eqnarray}
 |\nabla\vf_n(x)|&\le&\frac4{\sqrt{n}}\,|\eta'|_\9,\ \ff\,x\in\rrd,\label{e2.12}\\[1mm]
 |\Delta\vf_n(x)|&\le&\frac1n\,(2d|\eta'|_\9+8|\eta''|_\9),\ \ff\,x\in\rrd.\label{e2.13}
 \end{eqnarray} 
 \n By \eqref{e2.31a}, we have 
 \begin{equation}\label{e2.14}
 \hspace*{-5mm}\barr{l}
 \vf_n(y_1-y_2)-\lbb
 \Delta(\vf_n(\beta(y_1)-\beta(y_2)))
 +\lbb\,{\rm div}(D\vf_n(b^*(y_1)-b^*(y_2)))\vsp
 \qquad=\lbb(\nabla\vf_n\cdot D)(b^*(y_1)-b^*(y_2))
 -\lbb(\beta(y_1)-\beta(y_2))\Delta\vf_n\vsp
 \qquad-\,2\lbb\nabla\vf_n\cdot\nabla(\beta(y_1)-\beta(y_2))
 \mbox{ in }\cald'(\rr^d).\earr
 \end{equation}
Let $\calx_\delta:\rr\to\rr$ be the function
 \begin{equation}\label{e2.40a}
 	\calx_\delta(r)=\left\{\barr{rl}
  \dd\frac r\delta&\mbox{ if }|r|\le\delta,\vsp
  1&\mbox{ if }r>\delta,\vsp
  -1&\mbox{ if }r<-\delta,\earr\right.\end{equation}and let
  $$j_\delta(s)=\int^s_0\calx_\delta(r)dr,\ \ff\,s\in\rr.$$
 We know by \eqref{e2.4}--\eqref{e2.7prim} that
 $$\beta(y_i)\in L^p_{\rm loc},\ \ff p\in\mbox{$\left[1,\frac d{d-2}\)$},\
 \nabla\beta(y_i)\in L^q_{\rm loc},\
 \ff q\in\mbox{$\left[1,\frac d{d-1}\)$}.$$
 Moreover, by Hypothesis (iv) it follows that
 $$ |b^*(y_i)|\le\alpha_2|\beta(y_i)|,\ \ 
  |\nabla b^*(y_i)|\le \alpha_2|\nabla\beta(y_i)|,\ \ \mbox{ a.e. on }\rr^d,\ i=1,2,$$and, therefore, 
 $$b^*(y_i)\in L^p_{\rm loc},\ |\nabla b^*(y_i)|\in L^q_{\rm loc},\ i=1,2,\ \ff\,p\in\mbox{$\left[1,\frac d{d-2}\)$}, q\in\mbox{$\left[1,\frac d{d-1}\)$}.$$ This implies that 
$${\rm div}(Db^*(y_i))=D\cdot\nabla b^*(y_i)+b^*(y_i){\rm div}\,D\in L^1_{\rm loc},\ i=1,2,$$
because, by (ii), ${\rm div}(D)\in L^m_{\rm loc}$ for some $m>\frac d2$. Since $y_i\in D(A_0)$, $i=1,2,$ we have therefore  that
 $\Delta(\vf_{n}(\beta(y_1)-\beta(y_2)))$ and $\vf_n(\beta(y_1)-\beta(y_2))$ are in $L^1(\rrd)$. This yields
 $$\barr{l}
 \dd-\int_\rrd\Delta(\vf_n(\beta(y_1)-\beta(y_2)))
 \calx_\delta(\beta(y_1)-\beta(y_2))dx\\
 \qquad\dd=\int_{\rr^d}\nabla(\vf_n(\beta(y_1)-\beta(y_2)))\cdot\nabla(\beta(y_1)-\beta(y_2))\calx'_\delta(\beta(y_1)-\beta(y_2))dx\vsp
 \qquad\dd\ge\frac1\delta\int_{[|\beta(y_1)-\beta(y_2)|\le\delta]}
 (\nabla\vf_n\cdot\nabla(\beta(y_1)-\beta(y_2)))
 (\beta(y_1)-\beta(y_2))dx,\earr$$and so, by \eqref{e2.12}--\eqref{e2.14}, we have
 \begin{equation}\label{e2.15}
 \barr{l}
 \dd\int_\rrd\vf_n(y_1-y_2)\calx_\delta(\beta(y_1)-\beta(y_2))dx\\
 \qquad+\dd\frac\lbb\delta\!\int_{[|\beta(y_1)-\beta(y_2)|\le\delta]}
  \!\vf_n(b^*(y_1){-}b^*(y_2))
  (D\!\cdot\!\nabla(\beta(y_1){-}\beta(y_2)))dx\vsp 
   \qquad\le\dd\frac{C\lbb}{\sqrt{n}}
   +I^\delta_{\lbb,n},\earr\end{equation}where     $$I^\delta_{\lbb,n}\le\lbb|D|_\9 \int_{[|\beta(y_1)-\beta(y_2)|\le\delta]}|\vf_n|\cdot|\nabla(\beta(y_1)-\beta(y_2))|dx\to0\mbox{\ \  as\ \ }\delta\to0,$$because
    $$|\nabla \beta(y_1)-\nabla\beta(y_2)|=0,\mbox{ a.e. on $\{x;\,|\beta(y_1)(x)-\beta(y_2)(x)|=0\}.$}$$
       To obtain \eqref{e2.15},  we have used the relation
      $$\barr{l}
  -2\lbb\dd\int_{\rr^d}(\nabla\vf_n\cdot\nabla
  (\beta(y_1)-\beta(y_2)))\calx_\delta
  (\beta(y_1)-\beta(y_2))dx\\
    \qquad=-2\lbb\dd\int_{\rr^d}\nabla\vf_n\cdot\nabla j_\delta(\beta(y_1)-\beta(y_2)) dx\vsp
   \qquad=2\lbb\dd\int_{\rr^d}\Delta\vf_n j_\delta(\beta(y_1)-\beta(y_2))dx\vsp
  \qquad\le\dd\frac{2\lbb}{n}(2d|\eta'|_\9+8|\eta''|_\9)\dd\int_{\rr^d}(|\beta(y_1(x))|+|\beta(y_2(x))|)dx\vsp
  \qquad\le\dd\frac{4\lbb}{n}\,\alpha_1(d|\eta'|_\9+4|\eta''|_\9)
  (|y_1|_1+|y_2|_1)\le\dd\frac{C\lbb}{\sqrt{n}}\,|f|_1,\earr$$where $C$ is independent of $n$. 
  On the other hand, recalling that
 $$ D\in L^\9(\rr^d),\ \nabla \beta(y_i)\in L^q_{\rm loc},\  1\le q<\frac d{d-1},$$while by Hypothesis~(iv), we have
 $$|b^*(y_1)-b^*(y_2)|\le\alpha_2|\beta(y_1)-\beta(y_2)|,\mbox{ a.e. in }\rrd,$$and so, by \eqref{e2.15},  it follows for $\delta\to0$ that
 $$\int_\rrd\vf_n(y_1-y_2){\rm sign}(\beta(y_1)-\beta(y_2))dx\le\frac{C\lbb}{\sqrt{n}}\,|f|_1.$$
This yields  for $n\to\9$ that
 $$|y_1-y_2|_1=0,$$as claimed. $\Box$

 \begin{lemma}\label{l2.4} Assume that $d\ge1$. Then, for each $f\in L^1(\rrd)$ and all
 	$\lbb>0$, equation \eqref{e2.3} has a unique solution $y=J_\lbb(f)$. Moreover, one has
 	\begin{eqnarray}
 	&|J_\lbb(f_1)-J_\lbb(f_2)|_1\le|f_1-f_2|_1,\ \ff f_1,f_2\in L^1,\ \lbb>0, \label{e2.16}\\
 	&J_\lbb(f)\in\calp,\ \ff\,f\in\calp,\ \ff\lbb>0,\label{e2.17}
 	\end{eqnarray}
 	and $\ov{D(A_0)}=L^1$.
\end{lemma}

\n{\bf Proof.} The proof of Lemma \ref{l2.4} was given under the assumptions (i)--(iii) in \cite{4} (see also \cite{5}), so here its proof under our  assumptions (i)--(iii)  will be outlined only. 

We assume first that $f\in L^1\cap L^2$ and approximate equation \eqref{e2.3} by
\begin{equation}\label{e2.18}
y+\lbb(\vp I-\Delta)(\beta(y)+\vp y)+\lbb\,{\rm div}(D_\vp b^*_\vp(y))=f,\end{equation}where $b^*_\vp\in C^1_b(\rr)\cap C_b(\rr)$ is a smooth approximation of $b^*$ such that 
$$\mbox{$|b^*_\vp(r)|\le C|r|,$ $\lim\limits_{\vp\to0}b^*_\vp(r)=b(r)r$ uniformly on compacts,}$$and
$$D_\vp=\eta_\vp D,\ \eta_\vp\in C^1_0(\rr^d),\ 0\le\eta_\vp\le1,\ |\nabla \eta_\vp|\le1,\ \eta_\vp(x)=1\mbox{ if }|x|<\frac1\vp.$$Clearly, we have
\begin{equation}\label{e2.44}
\barr{c}
|D_\vp|\in L^\9\cap L^2,\ |D_\vp|\le|D|,\ \lim\limits_{\vp\to\9}D_\vp(x)=D(x),\mbox{ a.e. }x\in\rr^d.\vsp 
{\rm div}\ D_\vp\in L^1,\ ({\rm div}\,D_\vp)^-\le({\rm div}\,D)^-+\one_{\left[|x|>\frac1\vp\right]}|D|.\earr
\end{equation}
A typical example for $b _\vp$ is
$$\barr{c}
b_\vp\equiv b*\vf_\vp,\ \ \dd b^*_\vp(r)\equiv\frac{b_\vp(r)r}{1+\vp|r|},\ \ \ r\in\rr,\vsp
\dd\vf_\vp(r)=\frac1{\vp}\,\vf\(\frac r\vp\),\ \vf\in C^\9_0(\rr),\ \int_{\rr}\vf(x)dx=1.\earr$$
  \n We can rewrite \eqref{e2.18} equivalently as the following equation on $L^2$:
\begin{equation}\label{e2.19}
\lbb(\beta(y)+\vp y)+(\vp I-\Delta)\1 y+\lbb(\vp I-\Delta)\1{\rm div}(D_\vp b^*_\vp(y))=(\vp I-\Delta)\1f.\end{equation}We set
$$F(y)=\lbb(\beta(y)+\vp y)+(\vp I-\Delta)\1y+\lbb(\vp I-\Delta)\1{\rm div}(D_\vp b^*_\vp(y))$$ 
and note that   
$$\barr{l}
(F(y_1)-F(y_2),y_1-y_2)_2=\vp|(\vp I-\Delta)\1(y_1-y_2)|^2_2\vsp
\qquad +\lbb(\beta(y_1)-\beta(y_2),y_1-y_2)_2+\vp\lbb|y_1-y_2|^2_2\vsp
 \qquad+|\nabla(\vp I{-}\Delta)\1(y_1-y_2)|^2_2{-}\lbb(D_\vp (b^*_\vp(y_1){-}b^*_\vp(y_2)),\nabla(\vp I{-}\Delta)\1(y_1{-}y_2))_2\vsp
\qquad \ge\vp\lbb|y_1-y_2|^2_2+|(\vp I-\Delta)\1(y_1-y_2)|^2_2  +|\nabla(\vp I-\Delta)\1(y_1-y_2)|^2_2\vsp
 \qquad-\lbb|D_\vp |_\9|b^*_\vp|_{\rm Lip}|y_1-y_2|_2|\nabla(\vp I-\Delta)\1(y_1-y_2)|_2\ge0,\mbox{ for $0<\lbb<\lbb_\vp$.}
\earr$$ 

It is also clear that $(F(y),y)_2\ge\alpha_\vp\lbb|y|^2_2$ and so, for $0<\lbb<\lbb_\vp$, $F$ is monotone, continuous and coercive on $L^2(\rrd)$. Hence, it is surjective and so equation \eqref{e2.19} has a solution $y_\vp\in L^2(\rrd)$ and the latter is then true for all $\lbb>0$ (see Propositions 3.1 and 3.2 in \cite{2a}). By \eqref{e2.18}, it follows also that $y_\vp,\beta(y_\vp)\in H^1(\rrd)$.  

If $f\in L^1\cap L^\9$, we have
\begin{equation}
\label{e2.48}
|y_\vp|_\9\le(1+||D|+({\rm div}\,D)^-|^{\frac12}_\9)|f|_\9,\ \ 0<\lbb<\lbb_0<1.
\end{equation}Indeed, by \eqref{e2.18} we see that, for $M_\vp=|({\rm div}\,D_\vp)^-|^{\frac12}_\9|f|_\9$ and $\lbb<\lbb_0$,
$$\barr{l}
(y_\vp-|f|_\9-M_\vp)-\lbb\Delta(\beta_\vp(y_\vp)-\beta_\vp(|f|_\9+M_\vp))\vsp 
\qquad+\lbb \vp(\beta_\vp(u_\vp)-\beta_\vp(|f|_\9+M_\vp))+\lbb\,{\rm div}(D_\vp(b^*_\vp(u_\vp)-b^*_\vp(|f|_\9+M_\vp)))\vsp 
\qquad\le f-|f|_\9-M_\vp-\lbb b^*_\vp(M_\vp+|f|_\9){\rm div}\,D_\vp\le 0,\earr$$where $\beta_\vp(r)=\beta(r)+\vp r.$

 Multiplying the above equation by $\calx_\delta((y_\vp-(|f|_\9+M_\vp))^+)$ and integrating over $\rrd$, we get as above, for $\delta\to0$, $|(y_\vp-|f|_\9-M_\vp)^+|_1\le0$ and, therefore, by \eqref{e2.44} and since $b\ge0$,
$$y_\vp\le(1+||D|+({\rm div}\,D)^-|^{\frac12}_\9)|f|_\9,\mbox{ a.e. in}\rrd.$$Similarly, one gets that
$$y_\vp\ge-(1+||D|+({\rm div}\,D)^-|^{\frac12}_\9)|f|_\9,\mbox{ a.e. in }\rrd,$$and so \eqref{e2.48} follows.\newpage

Let us denote   the solution to \eqref{e2.18} by $y^\vp_\lbb(f)$ and define $\beta_\vp(r):=\beta(r)+\vp r,$ $r\in\rr.$ Then, we  multiply the equation
$$\barr{r}
(y^\vp_\lbb(f_1)-y^\vp_\lbb(f_2))
+\lbb(\vp I-\Delta)(\beta_\vp(y^\vp_\lbb(f_1))-\beta_\vp(y^\vp_\lbb(f_2)))\vspace*{2,5mm}\\
+\lbb\,{\rm div}(D_\vp(b^*_\vp(y^\vp_\lbb(f_1))-b^*_\vp(y^\vp_\lbb(f_2)))=f_1-f_2
\earr$$by $\calx_\delta(\beta_\vp(y^\vp_\lbb(f_1))-\beta_\vp(y^\vp_\lbb(f_2)))$ and integrate over $\rrd$. We set for $\lbb\in(0,\lbb_1)$ and $\delta>0$
$$E^\vp_{\lbb,\delta}=\{x\in\rr^d;\ |\beta_\vp(y^\vp_\lbb(f_1)(x))-\beta_\vp(y^\vp_\lbb(f_2)(x))|\le\delta\}.$$Since  $|\beta_\vp(r)-\beta_\vp(\bar r)|\ge\vp|r-\bar r|$, $r,\bar r\in\rr$, and $b^*_\vp$ is Lipschitz, we have
\begin{eqnarray}
&&\hspace*{-19mm}\dd\int_\rrd(y^\vp_\lbb(f_1)-y^\vp_\lbb(f_2))\calx_\delta(\beta_\vp(y^\vp_\lbb(f_1))-\beta_\vp(y^\vp_\lbb(f_2)))dx
\label{e2.26}\\
&&\hspace*{-19mm}\le|f_1-f_2|_1+\lbb\dd\int_\rrd(b^*_\vp(y^\vp_\lbb(f_1))-b^*_\vp(y^\vp_\lbb(f_2)))\nonumber\\[1mm]
&&\hspace*{-19mm}\qquad\quad D_\vp\cdot\nabla(\beta_\vp(y^\vp_\lbb(f_1)-\beta_\vp(y^\vp_\lbb(f_2)))\calx'_\delta(\beta_\vp(y^\vp_\lbb(f_1))-\beta_\vp(y^\vp_\lbb(f_2)))dx\hspace*{-16mm}\nonumber\\[2mm]
&&\hspace*{-19mm}\le|f_1{-}f_2|_1+
\dd\frac{C_\vp\lbb}{\delta}\!\dd\int_{E^\vp_{\lbb,\delta}}\!\!\!\! |\beta_\vp(y^\vp_\lbb(f_1)){-}\beta_\vp(y^\vp_\lbb(f_2))|
|\nabla(\beta_\vp(y^\vp_\lbb(f_1))){-}\beta_\vp(y^\vp_\lbb(f_2))||D_\vp|dx\hspace*{-16mm}\nonumber\vspace*{-7mm}\\ 
&&\hspace*{-19mm}\le |f_1{-}f_2|_1+C_\vp\lbb|D_\vp|_2\(\dd\int_{E^\vp_{\lbb,\delta}}\!\!\!
|\nabla(\beta_\vp(y^\vp_\lbb(f_1))-\beta_\vp(y^\vp_\lbb(f_2)))|^2dx\)^{\frac12}.\hspace*{-16mm}\nonumber
\end{eqnarray}
Then, letting $\delta\to0$ in \eqref{e2.26} and recalling that
$${\rm sign}(\beta_\vp(y^\vp_\lbb(f_1))-\beta_\vp(y^\vp_\lbb(f_2)))={\rm sign}(y^\vp_\lbb(f_1)-y^\vp_\lbb(f_2)),
\mbox{ a.e. in $\rrd$,}$$ we get by monotone convergence 
\begin{equation}
\label{e2.52a}
|y^\vp_\lbb(f_1)-y^\vp_\lbb(f_2)|_1\le|f_1-f_2|_1.
\end{equation}
In particular, for all $f\in L^1\cap L^\9$,
\begin{equation}
\label{e2.52aa}
\sup_{\lbb,\vp>0}|y^\vp_\lbb(f)|_1\le|f|_1.\end{equation}
Recall that by \eqref{e2.48} we have for some $C_D\in(0,\9)$
\begin{equation}
\label{e2.52aaa}
\sup_{^{\lbb\in(0,\lbb_0)}_{\vp>0}}|y^\vp_\lbb(t)|_\9\le C_D|f|_\9.\end{equation}
Hence, multiplying \eqref{e2.18} by $\beta(y_\vp)$ and integrating over $\rr^d$, we see that, for some $C\in(0,\9)$ and all $\lbb\in(0,\lbb_0)$, $\vp\in(0,1)$,\newpage
\begin{equation}
\label{e2.52aaaa}\lbb|\nabla\beta_\vp(y^\vp_\lbb(f))|^2_2\le C\(|f|_\9+\sup_{|r|\le C_D|f|_\9}|\beta(r)|\)|f|_1.
\end{equation}Now, fix $\lbb\in(0,\lbb_0)$. Set $y_\vp:= y^\vp_\lbb(f)$. Then, by \eqref{e2.52aa}--\eqref{e2.52aaaa}, $\{\beta_\vp(y_\vp)\}$ is bounded in $H^1$ and $\{y_\vp\}$ is bounded in $L^2$.

This implies that $\{\beta_\vp(y_\vp)\}$ is compact in $L^2_{\rm loc}$ and, therefore, along a subsequence $\{\vp\}\to0$, we have
$$\barr{rcll}
y_\vp&\to&y&\mbox{ weakly in }L^2,\vsp 
\beta(y_\vp)&\to&v&\mbox{strongly in }L^2_{\rm loc},\vsp 
\nabla\beta_\vp(y_\vp)&\to&\nabla v&\mbox{ weakly in }L^2.\earr$$
Since the map $y\to\beta(y)$ is maximal monotone in $L^2(\calo)$ for every bounded, open $\calo\subset\rr^d$, it follows that $v(x)=\beta(y(x))$, a.e. $x\in\rr^d$.
Moreover, as $\beta$ is, by Hypothesis (i), continuous and strictly monotone, it follows that 
\begin{equation}\label{e2.51prim}
y_\vp\to y\mbox{ a.e. on }\rr^d,\end{equation}and, therefore, we have
$$b_\vp(y_\vp)y_\vp \to b(y)y, \mbox{ a.e. in $\rrd$},\quad
b_\vp(y_\vp)y_\vp \to b(y)y \mbox{ weakly in $L^2$},$$selecting another subsequence $\{\vp\}\to0$, if necessary. Then, letting $\vp\to0$ in \eqref{e2.18}, we see that
\begin{equation}\label{e2.23}
y-\lbb\Delta\beta(y)+\lbb\,{\rm div}(Db(y)y)=f\mbox{ in }H\1(\rrd).\end{equation}
Moreover, we see that $\beta(y)\in H^1(\rrd)$  and \eqref{e2.52aa}--\eqref{e2.52aaaa} hold for $y,\beta(y)$ replacing $y^\vp_\lbb(f)$ and $\beta(y^\vp_\lbb(f))$, respectively. 
We denote this solution to equation \eqref{e2.23} by $J_\lbb(f)$. We shall now prove 
\begin{equation}\label{e2.25}
|J_\lbb(f_1)-J_\lbb(f_2)|_1\le|f_1-f_2|_1,\ \ff f_1,f_2\in L^1\cap L^\9,\ \lbb\in (0,\lbb_0),\end{equation}which implies, in particular, that $J_\lbb(f)\in L^1,$ $\ff f\in L^1\cap L^\9$. 

Let $f\in L^1$ be arbitrary but fixed and let $\{f_n\}\subset L^1\cap L^\9$ be such that  $f_n\to f$ in $L^1$ as $n\to\9$. We set $y_n=y(\lbb,f_n)$, that is,
\begin{equation}\label{e2.27}
y_n-\lbb\Delta\beta(y_n)+\lbb\,{\rm div}(Db(y_n)y_n)=f_n\mbox{ in }H\1(\rrd).\end{equation}Then, by \eqref{e2.52a}, we have 
$$|y_n-y_m|_1\le|f_n-f_m|_1,\ \ff n,m\in\nn,$$\newpage\n  and so there is $y=\lim\limits_{n\to\9} y_n$ in $L^1$. Since $A_0$ is closed on $L^1$, it follows by \eqref{e2.27} that $y\in D(A_0)$ and $y+\lbb A_0 y=f$ and so $R(I+\lbb A_0)=L^1$, as claimed. The fact that $y$ is  the unique solution to $y+\lbb A_0 y=f$ follows by Lemma \ref{l2.3}. Denoting this solution $y$ by $J_\lbb(f)$, we obtain by \eqref{e2.25} that \eqref{e2.16} holds. 

If $f\in\calp\cap L^\9$, it follows by \eqref{e2.18} that $y_\vp\in\calp$. By \eqref{e2.52aaa}, \eqref{e2.51prim}, it follows that $y_\vp\to J_\lbb(f)$ in $L^1_{\rm loc}$ as $\vp\to0$. Then, the argument from the proof of Lemma 3.3 in \cite{4} implies that $y_\vp\to J_\lbb(f)$ in $L^1$, so $J_\lbb(f)\in\calp$. Finally, by density, \eqref{e2.17} follows.

To prove that $\ov{D(A_0)}=L^1$, it suffices to note that, by (ii), (iii), $C^\9_0(\rr^d)\subset D(A_0)$ (because ${\rm div}\,D\in L^1_{\rm loc}$ and $\beta\in C^2(\rr^d)$, $b\in C^1\cap C_b$).   This completes the proof of Lemma~\ref{l2.4}. $\Box$

\bk\n{\bf Proof of Theorem \ref{t2.1} (continued).} By Lemmas \ref{l2.3} and \ref{l2.4}, it follows that $A_0$ is $m$-accretive in $L^1(\rrd)$ and $\ov{D(A_0)}=L^1.$
Then, as mentioned earlier, the existence and uniqueness of a mild solution $\rho$ to \eqref{e1.1} follows by the Crandall \& Liggett theorem (see \cite{1}, p.~154). Moreover, by \eqref{e2.17} and \eqref{e1.17} it follows that $\rho(t)\in\calp$, $\ff t\ge0$, if $\rho_0\in\calp$.

We shall show now that $\rho$ is a distributional solution to \eqref{e1.1}. Since $\rho_h\to\rho$ in $L^1((0,T)\times\rr^d)$ as $h\to0$, we have along a subsequence $\{h\}\to0$
$$\beta(\rho_h)\to\beta(\rho),\ \ b^*(\rho_h)\to b^*(\rho),\mbox{ a.e. in }(0,\9)\times\rr^d.$$
Then, taking into account that $|\beta(\rho_h)|\le \alpha_1|\rho_h|$, a.e. in $(0,T)\times\rr^d$, it follows by a standard argument and by Hypothesis (iv) that
$$\beta(\rho_h)\to\beta(\rho),\ \ b^*(\rho)\to b^*(\rho)\mbox{ in }L^1((0,\9)\times\rr^d)$$as $h\to0$.  By \eqref{e1.10}, we have 
\begin{equation}\label{e2.57}
\barr{l}
\dd\int^\9_h\!\!\!\int_{\rr^d}\!
\Big(\frac1h(\rho_h(t,x)-\rho_h(t-h,x))\vf(t,x)-\beta(\rho_h(t,x))\\\qquad\cdot\Delta\vf(t,x)-(D(x)\cdot\nabla\vf(t,x))b^*(\rho_h(t,x))\Big)dx\,dt=0\vsp 
\qquad\hfill \ff\vf\in C^\9_0([0,\9)\times\rr^d).\earr\end{equation}Taking into account that
$$\barr{ll}
\dd\int^\9_h\!\!\!\int_{\rr^d}\!
\rho_h(t-h,x)\vf(t,x)dx\,dt
\!\!\!&=\dd\int^\9_0\!\!\!\int_{\rr^d}\!
\rho_h(t,x)\vf(t+h,x)dx\,dt\vsp 
&+\dd\int^h_0\!\!\!\int_{\rr^d}\!\rho_0(x)\vf(t+h,x)dx\,dt,
\earr$$
and, letting $h\to0$ in \eqref{e2.57}, it follows  that \eqref{e1.22} holds, as claimed. 

If $\rho_0\in L^1\cap L^\9$, it follows as in Theorem 2.2 in \cite{4} that $\rho\in L^\9((0,T)\times\rr^d)$, $\ff T>0$.
This completes the proof. $\Box$

\begin{remark}	\label{r2.6}\rm 	Analyzing the proof of Lemma \ref{l2.4}, one sees that, as far as concerns the existence of a solution to equation \eqref{e2.3}, Hypothesis \eqref{e1.5} can be dispensed with. Moreover, the condition $b\ge0$ in Hypothesis (iii), which was used in Lemma \ref{l2.4} to prove \eqref{e2.48}, is no longer necessary for existence of a solution to \eqref{e2.3} if \eqref{e1.5} is strengthen to $\gamma_1|r|\le|\beta(r)|\le\alpha_1|r|,$ $\ff r\in\rr.$   \end{remark}

\n{\bf Proof of Theorem \ref{t2.2}.}\sk

\n   Assume that Hypotheses (i)--(v) hold and let $\rho_0\in L^1\cap L^2$. We set  $g(r)\equiv\int^r_0\beta(s)ds$.  Then, multiplying \eqref{e1.8} by $\beta(\rho^{j+1}_h)\in H^1(\rrd)$ and integrating over $\rrd$, we get
\begin{equation}\label{e2.37}\barr{l}
	\dd\int_\rrd g(\rho^{j+1}_h(x))dx+h\int_\rrd|\nabla\beta(\rho^{j+1}_h(x))|^2dx\vsp
	\quad\le h\dd\int_\rrd b^*(\rho^{j+1}_h(x))D(x)\cdot\nabla\beta(\rho^{j+1}_h(x))dx+\int_\rrd g(\rho^{j}_h(x))dx,\vsp\hfill j=0,1,...,N.\earr\end{equation}
Next, by \eqref{e2.37} we get
$$\barr{l}
\dd\int_\rrd g(\rho_h(t,x))dx+\int^t_0\int_\rrd|\nabla\beta(\rho_h(s,x))|^2dx\,ds\vsp
\qquad\le\dd\int_\rrd g(\rho_0(x))dx+C\int^t_0\int_\rrd\rho^2_h(s,x)ds\,dx,\  \ff\,t\in(0,T),\earr$$and, since by assumption (v) we have
$\alpha_0r^2\le g(r)\le\alpha_1 r^2,$ $\ff r\in\rr$, we get
\begin{equation}
\label{e2.39}
|\rho_h(t)|^2_2+|\beta(\rho_h)|^2_{L^2(0,T;H^1)}\le C_T|\rho_0|^2_2,\ \ \ff\,h,\ t\in[0,T].\end{equation}
For $h\to0$,  $\rho_h\to\rho$ in $L^\9(0,T;L^1)$ and so, along a subsequence $\rho_h\to\rho$, a.e. in $(0,T)\times\rrd$, we have
\begin{equation}\label{e2.58a}
\beta(\rho_h)\to\beta(\rho),\mbox{ a.e. on }(0,T)\times\rrd.\end{equation}Then, by \eqref{e1.5}   it follows that $\{\beta(\rho_h)\}$ is Lebesgue equi-integrable on \mbox{$(0,T)\!\times\!\rr^d$.} Hence, by the generalized Lebesgue convergence theorem (see, e.g., \cite[Theorem 21.4]{8prim}) it follows that $\beta(\rho_h)\to\beta(\rho)$ in $L^1((0,T)\times\rr^d)$.  Then, we have along a subsequence $\{h\}\to0$
 
\begin{eqnarray}
\rho_h&\to&\rho\qquad\quad \, \mbox{strongly in }L^1((0,T)\times\rrd),\label{e2.61b}\\
&&\qquad\ \ \ \  \mbox{ weakly$^*$ in }L^\9(0,T;L^2),\\[2mm]
\beta(\rho_h)&\to&\beta(\rho)\quad \ \ \mbox{strongly in }L^1_{\rm loc}((0,T)\times\rrd),\label{e2.61a}\\
&&\qquad\ \ \ \  \mbox{ weakly$^*$ in }L^\9(0,T;L^2),\nonumber\\[2mm]
\nabla\beta(\rho_h)&\to&\nabla\beta(\rho)\ \ \ \mbox{weakly in }L^2(0,T;(L^2(\rrd))^d).\label{e2.61aa}
\end{eqnarray}Moreover, letting $h\to0$ in \eqref{e2.39}, we get
$$|\rho(t)|^2_2+\|\beta(\rho)\|^2_{L^2(0,T;H^1)}\le C_T|\rho_0|^2_2.$$
This yields
\begin{equation}\label{e2.61aaa}
\Delta\beta(\rho),{\rm div}(Db^*(\rho))\in L^2(0,T;H\1(\rrd)),\ \ff\,T>0.\end{equation}
By \eqref{e1.10}, we get
$$\barr{l}
\dd\int^T_h
\(\psi(t),\(\dd\frac{\rho_h(t)-\rho_h(t-h)}h\)\)_2dt
+\dd\int^T_h{}{\raise-1,2mm\hbox{$_{H^1}$}}\!\<\psi(t),A_0(\rho_h(t))\>_{H^{-1}}dt\vsp
\quad-\dd\int^h_0(\psi(t),\rho_h(t))_2dt+\int^T_{T-h}(\psi(t),\rho_h(t))_2dt=0,\
\ff\psi\in C^\9_0((0,T)\times\rrd).\earr$$Equivalently,
$$\barr{r}
\dd-\int^{T-h}_0 \(\frac{\psi(t+h)-\psi(t)}h\,,\rho_h(t)\)_2dt
+\int^T_h{}{\raise-1,2mm\hbox{$_{H^1}$}}\!\<\psi(t)A_0(\rho_h(t))\>_{H\1}dt\vsp
-\dd\int^h_0(\psi(t),\rho_h(t))_2+\int^T_{T-h}(\psi(t),\rho_h(t))dt=0.\earr$$Taking into account \eqref{e2.61b}, \eqref{e2.61aa}, we get for $h\to0$
$$-\int^T_0(\rho(t),\psi'(t))_2dt+
\int^T_0{}{\raise-1,2mm\hbox{$_{H^1}$}}\!\<\psi(t),A_0(\rho(t))\>_{H\1}dt=0,\ \ff\,\psi\in C^\9_0((0,T)\times\rrd)$$and, by \eqref{e2.61aaa},  $A_0(\rho)\in L^2(0,T;H\1(\rrd))$. This means that
$$\frac{d\rho}{dt}=\Delta\beta(\rho)-{\rm div}(Db^*(\rho))=-A_0(\rho)\mbox{\ \ in }L^2(0,T;H\1(\rrd)).$$Hence 
$$\frac{d\rho}{dt}\in L^2(0,T;H\1(\rrd)),$$and, therefore, $\rho:[0,T]\to H\1(\rrd)$ is absolutely continuous and satisfies \eqref{e2.1b}--\eqref{e2.3b}. This completes the proof. $\Box$\bk

\section{The uniqueness of   distributional  solutions to NFPE}\label{s3}
\setcounter{equation}{0}

In this section, we  shall prove   the uniqueness of distributional solutions to \eqref{e1.1} under the following Hypotheses: 
\bit\item[\rm(j)] $\beta\in C^1(\rr),\ \beta'(r)\ge0,\ \ff\,r\in\rr,\ \beta(0)=0.$
\item[\rm(jj)] $D\in L^\9(\rrd;\rrd).$
\item[\rm(jjj)] $b\in C^1(\rr).$
\item[\rm(jv)] For each compact $K\subset \rr$ there exists $\alpha_K\in(0,\9)$ such that $$|b'(r)r+b(r)|\le\alpha_K|\beta'(r)|,\ \ff\,r\in K,$$ or, equivalently,
$$|b(r)r-b(\bar r)\bar r|\le\alpha_K|\beta(r)-\beta(\bar r)|,\ \ff\,r,\bar r\in K.$$
\eit
We note that Hypotheses (j)--(jjj) are    weaker than (i)--(iii).

\begin{remark}\label{r3.0} \rm If $\beta'(r)>0,\ \ff\,r\in\rr,$ then (jv) is automatically fulfilled due to (jjj).
\end{remark}

\begin{theorem}\label{t3.1} Let $d\ge1,\,T>0$,  and let $y_1,y_2\in L^{\red{\9}}((0,T){\times}\rrd)$ be two distributional solutions to \eqref{e1.1} on $(0,T)\times\rr^d$ $($in the sense of \eqref{e1.22}$)$ such that  $y_1-y_2\in L^\9(0,T;H\1)$ and  
\begin{equation}
\label{e3.1b}\lim_{t\to0}\ {\rm ess} \sup_{\hspace*{-4mm}s\in(0,t)}|(y_1(s)-y_2(s),\vf)_2|=0,\ \ff\vf\in C^\9_0(\rrd).
\end{equation}	Then $y_1\equiv y_2$. 
\end{theorem} 

\n{\bf Proof.} Replacing, if necessary, the functions $\beta$ and $b$ by
$$\beta_N(r)=\left\{\barr{ll}
\beta(r)&\mbox{ if }|r|\le N,\vsp
\beta'(N)(r-N)+\beta(N)&\mbox{ if }r>N,\vsp
\beta'(-N)(r+N)+\beta(-N)&\mbox{ if }r<-N,\earr\right.$$ and 
$$ b_N(r)=\left\{\barr{ll}
 b(r)&\mbox{ if }|r|\le N,\vsp
 b'(N)(r-N)+b(N)&\mbox{ if }r>N,\vsp
 b'(-N)(r+N)+b(-N)&\mbox{ if }r<-N,\earr\right.$$where $N\ge\max\{|y_1|_\9,|y_2|_\9\}$, we may assume that 
\begin{equation}\label{e4.1}
\beta',b'\in C_b(\rr),\end{equation}and, therefore, by (j) and (jv) we have
\begin{eqnarray}
|b(r)r-b(\bar r)\bar r|&\le&\alpha_3|\beta(r)-\beta(\bar r)|,\ \ \ff\,r,\bar r\in\rr,\label{e4.2a}\\[1mm]
(\beta(r)-\beta(\bar r))(r-\bar r)&\ge&\alpha_4|\beta(r)-\beta(\bar r)|^2,\ \ \ff\,r,\bar r\in\rr,\label{e4.2}
\end{eqnarray}where $\alpha_3=\alpha_{[-N,N]}$, $\alpha^{-1}_4=|\beta'|_{L^\9(-N,N)}$. We set
\begin{equation}\label{e3.2c}
\barr{c}
\Phi_\vp(y)=(\vp I-\Delta)\1y,\ \ff\,y\in L^2,\vsp
z=y_1-y_2,\ w=\beta(y_1)-\beta(y_2),\ b^*(y_i)\equiv b(y_i)y_i,\ i=1,2.\earr\end{equation}
It is well known that $\Phi_\vp:L^p(\rr^d)\to L^p(\rr^d)$, $\ff p\in[1,\9]$ and
\begin{equation}\label{e3.3a}
\vp|\Phi_\vp(y)|_p\le |y|_p,\ \ \ff y\in L^p,\ \vp>0.
\end{equation}
Moreover, $\Phi_\vp:L^2\to H^2,\ \ff\vp>0$.  
 
By   \eqref{e1.22}, we have 
$$z_t-\Delta w+{\rm div}\,  D(b^*(y_1)-b^*(y_2))=0\mbox{ in }\cald'((0,T)\times\rrd).$$
We set 
\begin{equation}
\label{e3.6a}
z_\vp=z*\theta_\vp,\ w_\vp=w*\theta_\vp,\ \zeta_\vp=(D(b^*(y_1)-b^*(y_2)))*\theta_\vp,\end{equation}where $\theta\in C^\9_0(\rr^d),$ $\theta_\vp(x)\equiv\vp^{-d}\theta\(\frac x\vp\)$ is a standard mollifier. We note that $z_\vp,w_\vp,\zeta_\vp,\Delta w_\vp,{\rm div}\,\zeta_\vp\in L^2(0,T;L^2)$ and  we have
\begin{equation}
\label{e3.5prim}
(z_\vp)_t-\Delta w_\vp+{\rm div}\,\zeta_\vp=0\mbox{ in }\cald'(0,T;L^2) .\end{equation}
This yields  $\Phi_\vp(z_\vp),\Phi_\vp(w_\vp)\in L^2(0,T;H^2),{\rm div}\,\Phi_\vp(\zeta_\vp)\in L^2(0,T;L^2)$ and 
\begin{equation}\label{e3.2cc}\barr{ll}
(\Phi_\vp(z_\vp))_t=\Delta\Phi_\vp(w_\vp)-{\rm div}\Phi_\vp(\zeta_\vp)=\vp\Phi_\vp(w_\vp)-w_\vp-{\rm div}\Phi_\vp(\zeta_\vp)\vsp\hfill\mbox{ in }\cald'(0,T;L^2).\earr\end{equation}
By \eqref{e3.5prim}, \eqref{e3.2cc} it follows that $(z_\vp)_t=\frac d{dt}\,z_\vp$,   $(\Phi_\vp(z_\vp))_t=\frac d{dt}\,\Phi_\vp(z_\vp)\in L^2(0,T;L^2)$, where $\frac d{dt}$ is taken in the sense of $L^2$-valued vectorial distributions on $(0,T)$. This implies that $z_\vp, \Phi_\vp(z_\vp)\in H^1(0,T;L^2)$ and both $[0,T]\ni t\mapsto z_\vp(t)\in L^2$ and $[0,T]\ni t\to\Phi_\vp(z_\vp(t))\in L^2$ are absolutely continuous. (See, e.g., \cite{1}, p.~23.) As a matter of fact, as seen above, in particular we have
\begin{equation}\label{e3.5a}
\Phi_\vp(z_\vp)\in L^2(0,T;H^2)\cap H^1(0,T;L^2), \end{equation}
 and, since the spaces $H^2,L^2$ are in duality with the pivot space $H^1=(H^2,L^2)_{\frac12}$, modifying the function $t\to\Phi_\vp(z_\vp(t))$ on a subset of measure zero, we have
	\begin{equation}\label{e3.9b}
	\Phi_\vp(z_\vp)\in C([0,T];H^1),
	\end{equation}(see, e.g., \cite{1}, p. 25).  
We set 
$$h_\vp(t)=(\Phi_\vp(z_\vp(t)),z_\vp(t))_2,\ t\in[0,T],$$ and get, therefore,
\begin{eqnarray}
h'_\vp(t)\!&\!\!=\!\!&\!2(z_\vp(t),(\Phi_\vp(z_\vp(t)))_t)_2\label{e4.3}\\
\!&\!\!=\!\!&\!2(\vp\Phi_\vp(w_\vp(t)){-}w_\vp(t){-}{\rm div}\Phi_\vp(\zeta_\vp(t)),z_\vp(t))_2\nonumber\\
\!&\!\!=\!\!&\!2\vp(\Phi_\vp(z_\vp(t)),w_\vp(t))_2{+}2(\nabla\Phi_\vp(z_\vp(t)),\zeta_\vp(t))_2\nonumber\\
&&-2(z_\vp(t),w_\vp(t))_2,\mbox{ a.e. }t\in(0,T).\nonumber
\end{eqnarray} 
By \eqref{e3.2cc}--\eqref{e4.3} it also follows that \mbox{$t\to h_\vp(t)$} is absolutely continuous on $[0,T]$. 
Since, by \eqref{e4.2}, \eqref{e3.2c}, \eqref{e3.6a},
\begin{equation}\label{e3.9a}
(z_\vp(t),w_\vp(t))_2\ge\alpha_4|w_\vp(t)|^2_2+\gamma_\vp(t),
\end{equation}where
\begin{equation}
\label{e3.9aa}
\gamma_\vp(t):=(z_\vp(t),w_\vp(t))_2-(z(t),w(t))_2,\end{equation}we get, therefore, by \eqref{e4.3}    that 
\begin{equation}\label{e4.4}
\hspace*{-6mm}\barr{ll}
0\le h_\vp(t)\!\!
\le h_\vp(0+){+}2\vp\!\dd\int^t_0\!\!(\Phi_\vp(z_\vp(s)),w_\vp(s))_2ds
{-}2\alpha_4\!\dd\int^t_0
|w_\vp(s)|^2_2ds\\
+2\alpha_3|D|_\9\dd\int^t_0\!|\nabla\Phi_\vp(z_\vp(s))|_2|w_\vp(s)|_2
ds+2
\int^t_0|\gamma_\vp(s)|ds,\, \ff t\in[0,T].\earr\hspace*{-6mm}\end{equation}
Taking into account that 
 $t\to\Phi_\vp(z_\vp(t))$ has an $H^1$-conti\-nuous version on $[0,T]$ (which we shall consider from now on), there exists $f\in H^1$ such that 
 $$\lim_{t\to0}\Phi_\vp(z_\vp(t))=f\mbox{\ \ in }H^1.$$ 
 Furthermore, for every $\vf\in C^\9_0(\rrd)$, $s\in(0,T),$
 $$0\le h_\vp(s)\le|\Phi_\vp(z_\vp(s)){-}f|_{H^1}|z_\vp(s)|_{H\1}{+}|f{-}\vf|_{H^1}|z_\vp(s)|_{H\1}{+}|(\vf*\theta_\vp,z(s))_2|.$$Hence, by \eqref{e3.1b},
 $$\barr{ll}
0\le h_\vp(0+)\!\!
 &=\dd\lim_{t\downarrow0}h_\vp(t)
 =\lim_{t\to0}\
 {\rm ess}\sup_{\hspace*{-4mm}s\in(0,t)}h_\vp(s)\vsp
 &\le\left(\dd\lim_{t\to0}|\Phi_\vp(z_\vp(t))-f|_{H^1}+|f-\vf|_{H^1}\right)|z_\vp|_{L^\9(0,T;H\1)}\vsp&
 +\dd\lim_{t\to0}\ {\rm ess}\sup_{\hspace*{-4mm}s\in(0,t)}|(\vf*\theta_\vp,z(s))_2|=|f-\vf|_{H^1}|z_\vp|_{L^\9(0,T;H\1)}.\earr$$Since $C^\9_0(\rrd)$ is dense in $H^1(\rrd)$, we find 
 \begin{equation}
 \label{e3.7aa}
 h_\vp(0+)=0.
 \end{equation} 
 On the other hand, taking into account that, for a.e. $t\in(0,T)$,
\begin{equation}\label{e4.4a}
\vp\Phi_\vp(z_\vp(t))-\Delta\Phi_\vp(z_\vp(t))=z_\vp(t),\end{equation}we get that
\begin{equation}\label{e4.5}
\hspace*{-5mm}\vp|\Phi_\vp(z_\vp(t))|^2_2{+}|\nabla\Phi_\vp(z_\vp(t))|^2_2{=}(z_\vp(t),\Phi_\vp(z_\vp(t)))_2{=}h_\vp(t),\mbox{ a.e. } t\in(0,T).\hspace*{-5mm}\end{equation}
Hence, $h_\vp$ has a continuous version on $[0,T]$, which we shall consider from now on, and \eqref{e4.5} holds for all $t\in[0,T]$.

By \eqref{e4.4}, \eqref{e3.7aa}, we get, for every $t\in[0,T]$, 
$$\barr{ll}
0\le h_\vp(t)\!\!\!&\le\vp\dd\int^t_0
|\Phi_\vp(z_\vp(s))|^2_2ds
-(\alpha_4-\vp)\int^t_0|w_\vp(s)|^2_2ds\vsp 
&+\dd\frac{\alpha^2_3|D|^2_\9}{\alpha_4}\int^t_0|\nabla\Phi_\vp(z_\vp(s))|^2_2ds+2\int^t_0|\gamma_\vp(s)|ds,\earr$$and so, by \eqref{e4.5}, this yields, for $\vp\in(0,\alpha_4]$,
\begin{equation}\label{e3.18}
0\le  h_\vp(t)\le
\max\(1,\frac{\alpha^2_3|D|^2_\9}{\alpha_4}\)\!\int^t_0\!\! h_\vp(s)ds{+}2\!\int^t_0|\!\gamma_\vp(s)|ds,\, \ff t\in[0,T].\end{equation}
On the other hand, we see by \eqref{e3.9aa}  that $\gamma_\vp\to0$ in $L^1(0,T)$ as $\vp\to0$. Then, \eqref{e3.18} yields
\begin{equation}
\label{e3.19}
0\le h_\vp(t)\le\eta_\vp(t)\exp\(\max\(1,\frac{\alpha^2_3|D|^2_\9}{\alpha_4}\),t\),\ \ff t\in[0,T],
\end{equation}
where $\dd\lim_{\vp\to0}\eta_\vp(T)=0$ and hence $\lim\limits_{\vp\to0} h_\vp=0$ in $C([0,T]).$ 

By \eqref{e4.5}, \eqref{e3.19},  for every  $t\in[0,T]$, it follows that, as $\vp\to0$, $\nabla\Phi_\vp(z_\vp(t))$ $\to0$ in $L^2$ and $\vp|\Phi_\vp(z_\vp(t))|^2_2\to0$, hence the left hand side of \eqref{e4.4a} converges to zero in $\cald'(\rrd)$ for every $t\in[0,T]$ and, so $0=\lim\limits_{\vp\to0} z_\vp(t)=z(t)$ in $\cald'(\rrd)$ for a.e. $t\in[0,T]$, which implies $y_1\equiv y_2$.~$\Box$
\begin{remark}\label{r3.3a} \rm If hypotheses (j), (jv) are strenghten to \eqref{e4.2a}--\eqref{e4.2}, then, as easily follows by the proof, one can replace in Theorem \ref{t3.1} condition $y_1,y_2\in L^\9((0,T)\times\rrd)$ by $y_1,y_2\in L^2((0,T)\times \rrd)$.  
\end{remark}
\begin{remark}\label{r3.2prim}\rm Let $\rho\in L^1((0,T)\times\rrd)$ such that $\beta(\rho)\in L^1((0,T)\times\rrd)$ and $\rho$ is a  solution to \eqref{e1.22}. Then, it is elementary to check that
	$$\int_\rrd\rho(t,x)dx=\int_\rrd\rho_0(dx)\mbox{\ \ for }dt\mbox{-a.e. }t\in(0,T).$$	
	Hence, if $\rho_0$ is nonnegative and $\rho\ge0$, a.e. on $(0,T)\times\rrd$, it follows by Lemma 2.3 in \cite{23'}   that there exists a $dt\otimes dx$-version $\wt\rho$ of $\rho$ such that\break for $\wt\rho_t(dx):=\wt\rho(t,x)dx,$ $t>0$, and $\wt\rho_0(dx):=\rho_0(dx)$, the map $[0,T]\ni t\mapsto\wt\rho_t$ is narrowly continuous.\end{remark}
Then, Remark \ref{r3.2prim} implies the following consequence of Theorem \ref{t3.1}. 

\begin{corollary}\label{c3.2''} Let $\rho_0\!\in\!\calm(\rrd)$, $\rho_0$ nonnegative, and $y_1,y_2\!\in\! L^{\9}((0,T)\!\times\!\rrd)$ $\cap L^1((0,T)\times\rrd)$ be two nonnegative solutions to \eqref{e1.22}.  Then, $y_1\equiv y_2$.
\end{corollary}

\n{\bf Proof.} Let $\wt y_1,\wt y_2$ be the $dt\otimes dx$-versions from Remark \ref{r3.2prim}. Then, for every $\vf\in C^\9_0(\rrd)$,
$$\dd\lim_{t\to0}\ {\rm ess}\sup_{\hspace*{-4mm}s\in(0,t)}
|(y_1(s)-y_2(s),\vf)_2|
=\dd\lim_{t\to0}\left|\int_\rr (\wt y_1(s,x)-\wt y_2(s,x))\vf(x)dx\right|=0.$$So, \eqref{e3.1b} holds and Theorem \ref{t3.1} implies the assertion. $\Box$

\begin{remark}\label{r3.2'''}\rm Theorem 5.2 in \cite{4} asserts that, under the assumptions (k), (kk), (kkk) and \eqref{e5.3} on $\beta,b$ and $D$, specified in \cite[Section 4]{4}, there exists a distributional solution to \eqref{e1.22} above for every $\rho_0\in\calm(\rrd)$. If (in the notation of Remark \ref{r3.2prim}) $S(t)\rho_0:=\wt\rho_t$, $t\ge0$, denotes the narrowly continuous solution to \eqref{e1.22}, it was left as an open question in \cite{4} whether the semigroup property
\begin{equation}
\label{e3.30}
S(t+s)\rho_0=S(t)S(s)\rho_0;\ t,s>0,
\end{equation}holds (see Remark 5.3 in \cite{4}). If, in addition to the above mentioned\break assumptions from \cite{4} on $\beta,b$ and $D$, we assume that our condition (jv) holds, then Corollary \ref{c3.2''} implies that \eqref{e3.30} holds. Indeed, for $s>0$ fixed, both functions in $t\ge0$ from the left and right hand sides of \eqref{e3.30} are in\break \mbox{$L^\9((0,T)\times\rrd)\cap L^\9(0,T;L^1)$} by \cite[Theo\-rem 5.2, (5.6), (5.7)]{4} and \cite[Theo\-rem 2.2, (2.6)]{4}, respectively, and both solve \eqref{e1.22}. Hence, Corollary \ref{c3.2''} implies that they are equal.
\end{remark}

\section{Weak uniqeness of the corresponding\\ McKean--Vlasov SDE}\label{s.4}
\setcounter{equation}{0}
Our next aim is to prove weak uniqueness of the 
 McKean--Vlasov SDE \eqref{e1.4} corresponding to the \FP\ equation \eqref{e1.1} (resp. \eqref{e1.22}). Here, we allow the law of $X_0$ to be a general probability measure $\rho_0$, not necessarily absolutely continuous with respect to Lebesgue measure. By Corollary \ref{c3.2''}, it is immediate to prove that the time-marginal law densities of all weak solutions to \eqref{e1.4} coincide, provided they are bounded. To prove that also their laws on path space coincide, requires to also prove the so-called "linearized uniqueness". This was already noted in \cite{5a}, \cite{7b}. 
\begin{theorem}\label{t4.1} {\rm(Linearized Uniqueness)} Assume that Hypotheses {\rm(j)--(jv)} hold. Let $T>0$ and $u\in L^{\9}((0,T)\times\rrd)$. Let $y_1,y_2\in L^{2}((0,T)\times\rrd)$ such that $y_1-y_2\in L^\9(0,T;H\1)$ and $y_1,y_2$ are solutions to the following linearized version of \eqref{e1.22}	
	\begin{equation}\label{4.1}
	\barr{r}
\dd\int^\9_0\!\!\int_\rrd\!
\(\vf_t+\frac{\beta(u)}u\Delta\vf{+}b(u)D\!\cdot\!\nabla\vf\)\rho\,dx\,dt
+\!\dd\int_\rrd\!\!\vf(0,x)\rho_0(dx)=0,\vsp 
\ff\vf\in C^\9_0([0,T)\times\rrd),
\earr\end{equation}
	for some $\rho_0\in\calm(\rrd)$, where $\frac{\beta(0)}0:=\beta'(0)$, such that \eqref{e3.1b} holds. Then, $y_1\equiv y_2$.
	\end{theorem}

\n{\bf Proof.} First, we note that by (j)--(jv) we have:
\begin{eqnarray}
&\dd\frac{\beta(u)}u,b(u)\in L^\9((0,T)\times\rrd),\label{4.2}\\
&\dd|Db(u)|\le\alpha_7\,\frac{\beta(u)}u,\mbox{ a.e. on }(0,T)\times\rrd,\label{4.3}
\end{eqnarray}where $\alpha_7:=|D|_\9\alpha_{[-|u|_\9,|u|_\9]}.$ Now, we set
\begin{equation}
\label{4.4}
z:= y_1-y_2,\ \ 
w:=\dd\frac{\beta(u)}u\,(y_1-y_2). 
\end{equation}Then, we have, since $\frac{\beta(u)}u\ge0$,
\begin{equation}
\label{4.5}
w  z=\frac{\beta(u)}u\,|y_1-y_2|^2\ge\(\left|\frac {\beta(u)}u\right|_\9+1\)\1|w|^2,\mbox{ a.e.  on }(0,T)\times\rrd,\end{equation}and
\begin{equation}
\label{4.6}
|Db(u)z|\le\alpha_7|w|.\end{equation}
We set
\begin{equation}
\label{4.7}
z_\vp=z*\theta_\vp,\ w_\vp=w*\theta_\vp,\ \zeta_\vp=(Db(u)(y_1-y_2))*\theta_\vp,\end{equation}where $\theta_\vp$ is as in the proof of Theorem \ref{t3.1}. Now, by \eqref{4.2}--\eqref{4.6}, we can repeat the proof of the latter line by line for $z_\vp$, $w_\vp$ and $\zeta_\vp$ in \eqref{4.7} to obtain $y_1\equiv y_2$. $\Box$

\begin{corollary}
	\label{c4.2} Let $\rho_0\in\calm(\rrd)$, $\rho_0\ge0$, and $y_1,y_2\in (L^1\cap L^2)((0,T)\times\rrd)$ be two nonnegative solutions to \eqref{4.1}.   Then, $y_1\equiv y_2$.
\end{corollary} 

\n{\bf Proof.} The assertion follows from Theorem \ref{t4.1} by analogous arguments as in the proof of Corollary \ref{c3.2''}. $\Box$

\begin{theorem}
	\label{t4.3} Assume that Hypotheses {\rm(j)--(jv)} hold and let $T>0$. Let $X(t),\ t\ge0$, and $\wt X(t),\ t\ge0$, be two solutions to \eqref{e1.4} such that, for
	$$u(t,\cdot):=\frac{d\call_{X(t)}}{dx},\ \ 
	\wt u(t,\cdot):=\frac{d\call_{\wt X(t)}}{dx},
$$we have
\begin{equation}
\label{4.8}
u,\wt u\in L^{\9}((0,T)\times\rrd).
\end{equation}Then, $X$ and $\wt X$ have the same laws, i.e., $\mathbb{P}\circ X\1=\wt{\mathbb{P}}\circ\wt X\,\1.$\end{theorem}

\n{\bf Proof.} We first note that, by narrow continuity, $$u(0,x)dx=\rho_0(dx)=
(\mathbb{P}\circ X(0)\1)(dx)=
(\wt{\mathbb{P}}\circ\wt X(0)\1)(dx).$$ Furthermore, again by narrow continuity, \eqref{4.8} implies
that $u(t,\cdot),\wt u(t,\cdot)\in L^\9(\rrd)$, for all $t\in(0,T].$ Furthermore, by It\^o's formula, both $u$ and $\wt u$ satisfy the weak (nonlinear) \FP\ equation \eqref{e1.22}. Hence, since $u,\wt u\in L^\9(0,T;L^1)\subset L^1((0,T)\times\rrd)$, Corollary \ref{c3.2''} implies $u\equiv \wt u$. Furthermore, again by It\^o's formula, $\mathbb{P}\circ X\1$ and $\wt{\mathbb{P}}\circ\wt X\1$ satisfy the martingale problem with the initial condition $\rho_0$ for the linear Kolmogorov operator
$$L_u:=\frac{\beta(u)}u\Delta+b(u)D\cdot\nabla.$$ Hence, by Corollary \ref{c4.2}, the assertion follows from Lemma 2.12 in \cite{10}. 

Here, for $s\in[0,T]$, the set $\mathcal{R}_{[s,T]}$, which appears in that lemma, is chosen to be the set of all narrowly continuous, probability measure-valued solutions of \eqref{4.1} having, for each $t\in[s,T]$, $t>0$, a density $u(t,\cdot)\in L^\9(\rrd)$ such that $v\in L^\9((0,R)\times\rrd)$. $\Box$

\begin{remark}
\label{r4.5}\rm For existence of weak solutions to \eqref{e1.4} with time marginals in $L^\9((0,T)\times\rrd)$ given by a solution to \eqref{e1.22} for initial conditions in $L^1\cap L^\9$, we refer to \mbox{\cite[Theorem 6.1(a)]{4}} in the much improved recent arXiv version. In particular, the conditions stated there are weaker than conditions (i)--(iv) from Section~1 of the present paper. But, as seen,  under the latter stronger conditions we have even mild solutions $\rho$ to \eqref{e1.1} (and hence solutions to \eqref{e1.22}) 
 with initial conditions $\rho_0$ merely in $L^1$, as follows from Theorem \ref{t2.1}, and we also have $\beta(\rho)\in L^1((0,T)\times\rrd)$. Hence, by the general results in \cite[Section 2]{2}, we get a weak solution to \eqref{e1.4} with time-marginals $\rho(t,x)dx,$ $t\in[0,T]$. Furthermore, additionally assuming that (v) holds, by Theorem \ref{t2.2} we even have that this mild solution $\rho$ is a strong solution to \eqref{e1.1} in $H\1$ in the sense of \eqref{e2.3b} if $\rho_0\in L^1\cap L^2$.
\end{remark}
\begin{remark}\label{r4.5a}\rm If, instead of (j), (jv) we assume the stronger conditions \eqref{e4.2a}, \eqref{e4.2}, then it follows by Remark \ref{r3.3a} that Theorem \ref{t4.3} still holds if \eqref{4.8} is replaced by the following weaker condition:
\begin{equation}\label{e4.9} u,\wt u\in L^2((0,T)\times\rrd)\cap L^\9(0,T;H\1).
\end{equation}Indeed, since $u,\wt u$ are probability densities, we always have $\wt u,\wt u\in L^\9(0,T;L^1)$. So, \eqref{e4.9} is indeed weaker than \eqref{4.8}.
\end{remark}
\section{Weak differentiability of the nonlinear\\ \FP\ flow}\label{s4}
\setcounter{equation}{0}

Though the continuous semigroup $S(t):L^1\to L^1,$ $t\ge0$, defined by Theo\-rem \ref{t2.1} is not differentiable on $(0,\9)$, it is however differentiable in the distribution space $H\1=H\1(\rrd)$. Namely, consider in the space $H\1$ the operator
 $A_1:D(A_1)\subset H\1\to H\1$, 
$$\barr{rcl}
A_1(y)&{=}&-\Delta\beta(y)+{\rm div}(Db(y)y),\ \ \ff\,y\in D(A_1),\vsp
D(A_1)&{=}&\{y\in L^2;\beta(y)\in H^1(\rrd)\}.\earr$$
(We note that $D(A_1)$ is dense in $H\1$.)

We shall assume here that besides (i)--(iv) the following hypothesis holds.
\begin{itemize}
	\item[(vi)] $\beta\in {\rm Lip}(\rr).$
	\end{itemize}

\begin{lemma}\label{l4.1b} The operator $A_1$ is~quasi-$m$-accretive in $H\1$, that is, $\oo I+A_1$   is $m$-accretive in $H\1$ for some $\oo\ge0$.\end{lemma}

\n{\bf Proof.} One must prove that 
\begin{eqnarray}
&	\hspace*{-4mm}R(I+\lbb A_1)=H\1,\ \ \ff\,\lbb\in(0,\oo\1),\label{e4.1b}\\[1mm]
	&	\hspace*{-4mm}\|(I{+}\lbb A_1)\1u{-}(I{+}\lbb A_1)\1v\|_{H\1}\le(1{-}\lbb\oo)\1\|u{-}v\|_{H\1},\, \ff\,u,v\in H\1.\qquad\label{e4.2b}
	\end{eqnarray}
To prove \eqref{e4.1b}, consider for $f\in H\1$ the equation
\begin{equation}
	\label{e4.3b}
	y-\lbb\Delta\beta(y)+\lbb\,{\rm div}(Db(y)y)=f\mbox{ in }\cald'(\rrd),\end{equation}
and approximate it by \eqref{e2.18}. Arguing as in the proof of Lemma \ref{l2.3}, we get for the solution $y_\vp$ to \eqref{e2.18} the estimate 
$$|y_\vp|^2_2+\lbb|\nabla\beta(y_\vp)|^2_2|+\vp|\nabla y_\vp|^2_2\le C\|f\|^2_{H\1},\ \ \ff\vp>0,\ 0<\lbb<\lbb_0=\oo\1.$$
Then, letting $\vp\to0$, we get as above that there is the limit
$$\barr{rl}
y=\dd\lim\limits_{\vp\to0}y_\vp&\mbox{ strongly in $L^2_{\rm loc}$ and weakly in $L^2$,}\vsp
\beta(y_\vp)\to\beta(y)&\mbox{ strongly in $L^2_{\rm loc},$ weakly in $H^1$},\earr$$and
$$\barr{c}
y+\lbb A_1(y)=f,\vsp
 |y|^2_2+|\nabla\beta(y)|^2_2\le C\|f\|^2_{H\1}.\earr$$
Moreover, subtracting equations \eqref{e4.3b} for two solutions $y_1,y_2\in D(A_1)$ corres\-pon\-ding to $f_1,f_2\in H\1$, multiplying by $(I-\Delta)\1(y_1-y_2)$ and integrating on $\rrd$, we get
$$\barr{l}
\|y_1-y_2\|^2_{H\1}+\lbb(\beta(y_1)-\beta(y_2),y_1-y_2)_2\vsp
\qquad
\le\|f_1-f_2\|_{H\1}\|y_1-y_2\|_{H\1}+
\lbb\|\beta\|_{\rm Lip}|y_1-y_2|_2\|y_1-y_2\|_{H\1}\vsp
\qquad+\lbb|D|_\9|b^*(y_1)-b^*(y_2)|_2\|y_1-y_2\|_{H\1}.\earr$$This yields
$$\|y_1-y_2\|_{H\1}\le(1-\lbb\oo)\1\|f_1-f_2\|_{H\1},\ \ff\,\lbb\in(0,\oo\1),$$and so \eqref{e4.2b} follows.

By Lemma \ref{l4.1b}, we infer that the operator $-A_1$ generates a continuous quasi-contractive semigroup in $H\1$. More precisely, we have (see \cite{1}, p. 143)

\begin{theorem}\label{t4.1b} Assume that Hypotheses {\rm(i)--(iv)} hold. Then, for each $y_0\in D(A_1)$, there is a unique absolutely continuous function $y:[0,\9)\to H\1(\rrd)$ $t$-differentiable from the right, such that
\begin{equation}\label{e4.4b}
	\barr{l}
	\dd\frac d{dt}\ y(t)+A_1(y(t))=0,\ \ \mbox{ a.e. }t>0,\vsp
	y(0)=y_0,\earr
	\end{equation}
and 
\begin{equation}
	\label{e4.5b} \frac{d^+}{dt}\ y(t)+A_1(y(t))=0,\ \ \ff\,t\ge0.
	\end{equation}
Moreover, $S_1(t)y_0\equiv y(t)$ satisfies

\begin{eqnarray}
	&S_1(t+s)y_0=S_1(t)S_1(s)y_0,\ \ff\,t,s\ge0,\ y_0\in D(A_1),\label{e4.6b}\\[1mm]
	&\|S_1(t)y_0{-}S_1(t)\bar y_0\|_{H\1}\le\exp(\oo t)\|y_0{-}\bar y_0\|_{H\1},\, \ff t\ge0,\, y_0,\bar y_0\!\in\! D(A_1).\quad\label{e4.8b}\\[1mm]
	&S_1(t)y_0=\dd\lim_{n\to\9}\(I+\frac tn\,A_1\)^{-n}y_0\mbox{ in }H\1,\ \ff\,t\ge0.\label{e4.8c}
	\end{eqnarray}
\end{theorem}

The semigroup $\{S_1(t)\}_{t\ge0}$ extends by density on the closure $\ov{D(A_1)}=H\1$ of $D(A_1)$ in $H\1$ and by \eqref{e4.5b} it follows that $-A_1$ is the infinitesimal generator of $S_1(t)$ in $H\1$. Since $S(t)$ and $S_1(t)$ are given by the same exponential formula \eqref{e1.17} (respectively \eqref{e4.8c}) on $L^1$ (respectively $H\1$), we infer that

$$S(t)\rho_0=S_1(t)\rho_0,\ \ff\,t\ge0,\ \rho_0\in L^1\cap H\1.$$
This implies that the {\it\FP\ semigroup $S(t)$ given by Theorem {\rm\ref{t2.1}} is $H\1$  dif\-fe\-ren\-tiable from the right in $t$ on $L^1\cap H\1$.}

\section{Fokker--Planck equations with $x$-dependent diffusion coefficients}
\setcounter{equation}{0}

Theorems \ref{t2.1}, \ref{t2.2}, as well as Theorem \ref{t3.1}, extend {\it mutatis-mutandis} to nonlinear Fokker--Planck equations of the form  
\begin{equation}\label{e5.1}
\barr{l}
\rho_t(t,x)-\Delta\beta(x,\rho(t,x))+{\rm div}(D(x)b(\rho(t,x))\rho(t,x))=0,\\\hfill t\ge0,\ x\in\rr^d,\vsp
\rho(0,x)=\rho_0(x),\ x\in\rr^d,
\earr\end{equation}where $\beta:\rr^d\times\rr\to\rr,\ b:\rr\to\rr,\ D:\rr^d\to\rr$ satisfy (ii), (iii) and
\begin{itemize}
	\item[$\rm(i)'$] $\beta\in C^2(\rr^d\times\rr),\ \beta_r(x,r)>0,\ \ff r\ne0,\ x\in\rr^d,$ $\beta(x,0)\equiv0$ and
	\begin{eqnarray}
		&\label{e5.2}
	h_N(x)=\sup\{|\Delta_x\beta(x,r)|;\,|r|\le N\}\in L^1(\rr^d),\ \ff N>0,\\[2mm]
&|\beta(x,r)|\le\alpha_2|r|,\ \ff x\in\rr^d,\ r\in\rr.\label{e5.3}	\end{eqnarray}
		\item[$\rm(iv)'$]\ \vspace*{-9mm}
		\begin{equation}\label{e5.4}
					|b^*(r)-b^*(\bar r)|\le\alpha_2|\beta(x,r)-\beta(x,\bar r)|,\ \ff x\in\rr^d,\ r,\bar r\in\rr.
\end{equation}
\end{itemize} 
The definition of the mild solution $\rho$ to \eqref{e5.1} is that given in \eqref{e1.6}--\eqref{e1.9}. We have

\begin{theorem}\label{t5.1} Under Hypotheses $\rm(i)'$, {\rm(ii), (iii)}, $\rm(iv)'$, for each $\rho_0\in L^1$ there is a unique mild solution $\rho$ to \eqref{e5.1} which has all properties mentioned in Theorem {\rm\ref{t2.1}}.\end{theorem}
The proof is exactly the same as in the previous case and relies on Proposition \ref{p2.2}, where $A_0:L^1\to L^1$ is the operator
$$\barr{rcl}
A_0(y)&=&-\Delta\beta(x,y)+{\rm div}(D(x),b(y)y),\vsp
D(A_0)&=&\{y\in L^1;-\Delta\beta(x,y)+{\rm div}(D(x)b(y)y)\in L^1\}.\earr$$ 
The extension of Theorem \ref{t2.2}   to \eqref{e5.1} is also immediate, and so we omit the details. 
As regards the uniqueness of a distributional solution to equation \eqref{e5.1}, instead of (j) the following hypothesis will be assumed:
\begin{itemize}
	\item[$\rm(j)'$] $\beta,\beta_r\in C(\rr^d\times\rr), \beta_r\in L^\9(\rr^d\times(-N,N)),$ $\ff N>0,$ $\beta_r(x,r)\ge0,$ $\ff x\in\rr^d,$ $r\in\rr,$ $\beta(x,0)\equiv0.$
	\end{itemize}Then, we obtain results analogous to Theorem \ref{t3.1} and Corollary \ref{c3.2''}.

The proofs are identical   because also in this case one may assume that (see \eqref{e4.1})
$$\beta\in C_b(\rr^d\times\rr);\ b\in C_b(\rr),\ \beta_r\in L^\9(\rr^d\times\rr).$$We also note that Theorem \ref{t4.1b} remains valid in this case if, besides Hypotheses $\rm(i)'$, (ii), (iii), $\rm(iv)'$, one assume that
$$\beta_r\in L^\9(\rr^d\times\rr).$$Then, it follows as in Lemma \ref{l4.1b} that \eqref{e4.1b}--\eqref{e4.2b} hold, and so, the operator $A_1$ generates in $H\1(\rr^d)$ a semigroup $S_1(t)$ satisfying \eqref{e4.4b}--\eqref{e4.8c}.

As in Section \ref{s.4},  
these results imply weak uniqueness results for probabilistically weak solutions to the McKean--Vlasov equation
\begin{equation*}
\barr{l}
dX(t)=D(X(t))b(\rho(t,X(t)))dt+\dd \sqrt{\frac{2\beta(X(t),\rho(t,X(t)))}{\rho(t,X(t))}}\ dW(t),\vsp
X(0)=X_0.\earr\end{equation*}
Details on this will be given a forthcoming survey paper which is in preparation.

 \n{\bf Acknowledgements.} The authors are indebted  to Haim Brezis for useful suggestions regarding the proof of estimates \eqref{e2.7prim}.

\n This work was supported by the DFG through SFB 1283/2 2021-317210226 and by a grant of the Ministry of Research, Innovation and Digitization, CNCS--UEFISCDI project  PN-III-P4-PCE-2021-0006, within PNCDI III.

\end{document}